\documentclass{article}

\usepackage{PRIMEarxiv}

\usepackage[utf8]{inputenc} 
\usepackage[T1]{fontenc}    
\usepackage[colorlinks,linkcolor=blue,citecolor=red,urlcolor=blue]{hyperref}       
\usepackage{url}            
\usepackage{booktabs}       
\usepackage{amsfonts}       
\usepackage{nicefrac}       
\usepackage{microtype}      
\usepackage{lipsum}
\usepackage{fancyhdr}       
\usepackage{graphicx}       
\graphicspath{{media/}}     

\usepackage{amsmath}
\usepackage{amsfonts}
\usepackage{amssymb}
\usepackage{amsthm}
\usepackage[upright]{fourier} 
\usepackage{mathtools}
\usepackage{ifthen}
\usepackage{lettrine}
\usepackage[capitalise]{cleveref} 
\usepackage{tikz}
\usetikzlibrary{calc, math}
\usetikzlibrary{positioning}
\usepackage{tikz-cd}

\newtheorem{theorem1}{Theorem}[section]
\Crefname{theorem1}{Theorem}{Theorems}

\Crefname{lemma}{Lemma}{Lemmas}
\newtheorem{proposition1}{Proposition}[section]
\Crefname{proposition1}{Proposition}{Propositions}

\theoremstyle{definition}

\newtheorem{assumption1}{A}[section]
\Crefname{assumption1}{Assumption}{Assumptions}

\theoremstyle{remark}
\newtheorem{remark}{Remark}[section]

\numberwithin{equation}{section}

\newcommand{\R}{\ensuremath{\mathbb{R}}}

\definecolor{tcb}{rgb}{0.1,0.3,1}
\usepackage{orcidlink}

\usepackage{ccicons} 
\usepackage{academicons} 
\usepackage{fontawesome} 
\usepackage{pifont} 
\usepackage{orcidlink}
\usepackage[hyperpageref]{backref}  

\usepackage[numbers]{natbib} 

\newcommand{\mapolicebackref}[1]{
    \hspace*{\fill} \mbox{\textit {\small #1}}
}
\renewcommand*{\backref}[1]{}
\renewcommand*{\backrefalt}[4]{%
\ifcase #1 \mapolicebackref{no quotes}
    \or \mapolicebackref{quoted page #2}
    \else \mapolicebackref{#1 quotes pages #2}
\fi
}

\definecolor{fondtitre}{RGB}{85,85,85}
\definecolor{fonddeboite}{RGB}{232,232,232}
\definecolor{shadecolor}{RGB}{232,232,232}

\title{\Huge\textbf{Exponential ergodicity of mean-field Langevin dynamics by synchronous coupling}
\thanks{\textit{\underline{Acknowledgements}}: 
\textbf{The \href{https://mon-portfolio-de-chercheur.webnode.fr/}{author} thanks \href{https://blog.univ-angers.fr/panloup/}{Fabien \textsc{Panloup}} for extensive discussions and suggestions. This paper was originally released as a preprint on the HAL repository, and constitutes the second chapter of the \href{https://theses.hal.science/tel-05460067}{author’s Ph.D. thesis}. The author is grateful to the reviewers for their valuable comments during the revision process. Thanks also to everyone who has contributed in one way or another to the realization of this work. Last but not least, thankful for the benefits from \href{https://www.lebesgue.fr/en}{\textsc{Henri Lebesgue Center}} (program \textbf{ANR-11-LABX-0020-01}) such as \href{https://www.lebesgue.fr/en/content/bourses_master}{\textsc{Master Scholarships}}, \href{https://www.lebesgue.fr/en/node/4878}{\textsc{Lebesgue Doctoral Meeting}}, \href{https://www.lebesgue.fr/en/content/post-doc}{\textsc{Post doc positions}}... and thankful for the benefits for the ANR project \href{https://math.univ-angers.fr/~chaumont/rawabranch/index.html}{\textsc{Rawabranch}} number \textbf{ANR-23-CE40-0008}.}\newline 
\textit{\underline{Mathematics Subject Classification}}: \textbf{39B62; 82C31; 26D10; 47D07; 60G10; 60H10; 60J60}}
}

\author{
  \Huge{\orcidlink{0000-0002-3281-1954}} \Large{\href{https://mon-portfolio-de-chercheur.webnode.fr/}{Mohamed Alfaki ABOUBACRINE ASSADEK}} \\
   \href{http://www.univ-angers.fr/}{Univ Angers,} \href{https://www.cnrs.fr/fr}{CNRS,} \href{http://recherche.math.univ-angers.fr/}{LAREMA,}\\
   \href{https://sfrmathstic.univ-angers.fr/fr/index.html}{SFR MATHSTIC,}\\
   \href{https://www.angers.fr/}{F-49000 Angers, France }\\
  \texttt{\href{mailto:mohamedalfaki.agaboubacrineassadeck@univ-angers.fr}{mohamedalfaki.agaboubacrineassadeck}@univ-angers.fr} \\
  \AND
  \today\\
}

\begin{document}
\maketitle

\begin{abstract}
In the continuity of \cite{assadeck2023exponentialergodicitymckeanvlasovsde,BolleyGuillinMalrieu}, the paper investigates both the long-term quantitative behavior of first-order (overdamped) and second-order (underdamped) mean-field Langevin dynamics and their uniform-in-time propagation of chaos. First, sufficient conditions (involving small interaction) are provided regarding the macroscopic Hamiltonian and its flat derivatives to ensure the contraction of the invariant measure map for the overdamped mean-field Langevin process. Secondly, using synchronous coupling, we demonstrate exponential contraction in the Wasserstein $\mathcal{W}_2$ metric for both overdamped and underdamped mean-field Langevin dynamics. Next, we deduce uniform-in-time $L^2$-propagation of chaos with a Fournier–Guillin bound. In the kinetic case, general drift conditions are allowed
in the velocity component resulting in a non-explicit representation of a unique invariant measure. The results are applied to several concrete and original examples.
\end{abstract}

\vfill

\newpage
\tableofcontents

\addcontentsline{toc}{section}{Introduction}
\section*{Introduction}
In the continuity of \cite{assadeck2023exponentialergodicitymckeanvlasovsde,BolleyGuillinMalrieu}, the present work studies both the quantitative long-time behavior of the first order (over-damped) and second order (under-damped)
mean field Langevin dynamics and their uniform-in-time propagations of chaos properties, under a functional not necessary convex in the flat interpolation sense\footnote{convexity of $H$ in the flat interpolation sense :\begin{equation*}
    \forall (t,\mu,\nu)\in[0,1]\times\mathcal{P}(\mathbb{R}^{d})\times\mathcal{P}(\mathbb{R}^{d}),\quad H((1-t)\mu+t\nu)\leq (1-t)H(\mu)+tH(\nu).
\end{equation*}}: we consider the first order McKean-Vlasov equation
\begin{equation}\label{MLD}
    dX_{t}=\sigma dB_{t}-\mathcal{D}_{m}H(\mu_{t},X_{t})dt;
\end{equation}
and the second order McKean-Vlasov equation
\begin{equation}\label{KineticMLD}
    \begin{cases}
            dP_{t}&=V_{t}dt;\\
            dV_{t}&=\sigma dB_{t}-A(V_{t})dt-B(P_{t})dt-\mathcal{D}_{m}H(\mathbb{P}_{P_{t}},P_{t})dt;
        \end{cases}
\end{equation}
where $\mathcal{D}_{m}H$ denotes the intrinsic derivative of the functional $H$ (a mapping from $\mathcal{P}_{p}(\mathbb{R}^{d})\times\mathbb{R}^{d}$ into $\mathbb{R}^{d}$ and will defined in \cref{intrinsecderiv}), $\sigma>0$,  $\mu_{t}:=\mathbb{P}_{X_{t}}$, $A$ a friction field (velocity field) and $B$ an external confinement field. Without loss of generality, in the remainder of the paper, we will take $\sigma=\sqrt{2}$.\newline 

Mean-field models are the topic of active research. In particular, mean-field Langevin dynamics recently attracts attention due to its connection to (noisy) \textit{gradient descent} on infinitely wide \textit{neural networks} in the mean-field regime, and hence the convergence property of the dynamics is of great
theoretical interest: see, among others, \cite{aboubacrineassadek:tel-05460067,assadeck2023exponentialergodicitymckeanvlasovsde,chen2023uniformintimepropagationchaosmean,Chen_2024,monmarché2023logarithmicsobolevinequalitiesnonequilibrium,Monmarch__2024,hu2020meanfieldlangevindynamicsenergy,wang2024uniformlogsobolevinequalitiesmean,kazeykina2023ergodicityunderdampedmeanfieldlangevin} for more details.\newline

Let $n\geq1$ be an integer and let $\mathbf{x}:=(x_{1},\ldots,x_{n})\in(\mathbb{R}^{d})^{n}$  be an $n$-tuple of coordinates in $\mathbb{R}^{d}$. We denote by $\mu_{\mathbf{x}}$ the empirical measure formed with the $n$ coordinates $x_{1},\ldots,x_{n}$, that is,
\begin{equation}
    \mu_{\mathbf{x}}:=\frac{1}{n}\sum_{i=1}^{n}\delta_{x_{i}}.
\end{equation}
We consider the first order mean field particle system  
\begin{equation}
    \forall i\in\{1,\ldots,n\},\quad dX^{(i)}_{t}=\sigma dB^{(i)}_{t}-\mathcal{D}_{m}H(\mu_{\mathbf{X}_{t}},X^{(i)}_{t})dt;\\
\end{equation}
and the second order mean field particle system
\begin{equation}
    \begin{cases}
            dP^{(i)}_{t}&=V^{(i)}_{t}dt;\\
            dV^{(i)}_{t}&=\sigma dB^{(i)}_{t}-A(V^{(i)}_{t})dt-B(P^{(i)}_{t})dt-\mathcal{D}_{m}H(\mu_{\mathbf{P}_{t}},P^{(i)}_{t})dt;
        \end{cases}
\end{equation}
where $B^{(1)},\ldots,B^{(n)}$ are i.i.d. standard Brownian motions in $d$-dimensions. Under suitable assumptions, we have (\textit{propagation of chaos})
\begin{align}
    \frac{1}{n}\sum_{i=1}^{n}\delta_{X^{(i)}_{t}}&\overset{n\to+\infty}{\longrightarrow}\mu_{t},\\
    \frac{1}{n}\sum_{i=1}^{n}\delta_{(P^{(i)}_{t},V^{(i)}_{t})}&\overset{n\to+\infty}{\longrightarrow}\mathbb{P}_{(P_{t},V_{t})}
\end{align}
at least in the sense of distributions. We refer readers to \cite{reviewpropachaosI,reviewpropachaosII} for a more
comprehensive review of propagation of chaos. For $\mu\in\mathcal{P}(\mathbb{R}^{d})$, we consider the over-damped Langevin dynamics:
\begin{equation}
    dX^{\mu}_{t}=\sqrt{2}dB^{\mu}_{t}-\mathcal{D}_{m}H(\mu,X^{\mu}_{t})dt.
\end{equation}
By explicit computations, we can find that the Gibbs measure:
\begin{equation}
    \Phi(\mu)(dx):=\frac{e^{-\frac{\delta H}{\delta m}(\mu,x)}}{\int_{\mathbb{R}^{d}}e^{-\frac{\delta H}{\delta m}(\mu,y)}dy}dx
\end{equation}
is invariant to the dynamics and $\frac{\delta H}{\delta m}$ is called \textit{flat derivative} of $H$ and is a continuous mapping from $\mathcal{P}_{p}(\mathbb{R}^{d})\times\mathbb{R}^{d}$ into $\mathbb{R}$ and will defined in \cref{deffirstflatderivative}. By explicit computations, we can find that the set of invariant probability distributions of the mean field Langevin dynamics (over-damped) is defined by (\textit{fixed-point problem})
\begin{equation}
    \{\mu_{\infty},\quad\Phi(\mu_{\infty})=\mu_{\infty}\}.
\end{equation}
In the kinetic case, in general, if the \textit{velocity} and \textit{external confinement fields} do not derive from \textit{potential gradients}, the \textit{invariant probabilities} are not explicit.\newline

As for the \textit{long-time behavior}, the \textit{uniform in time propagation of chaos} and the \textit{convergence rates}, many works have been devoted to the subject and different approaches have been introduced such as \textit{coupling methods} and \textit{functional inequalities} (see, among others, \cite{aboubacrineassadek:tel-05460067,assadeck2023exponentialergodicitymckeanvlasovsde,chen2023uniformintimepropagationchaosmean,Chen_2024,BGL-MarkovDiffusion,BolleyGuillinMalrieu,monmarché2023logarithmicsobolevinequalitiesnonequilibrium,malrieu:hal-01282602,hu2020meanfieldlangevindynamicsenergy,Monmarch__2024,wang2024uniformlogsobolevinequalitiesmean,guillin2019kineticfokkerplanckequationmean,GuillinWuZhang,cattiaux2006probabilisticapproachgranularmedia}).\newline 

Let us briefly describe the \textit{coupling method} (the \textit{major approach of this paper}). The basic idea is that an upper bound on the Wasserstein
distance between two probability distributions is given by the construction of any pair of random variables distributed respectively according to those. The goal is thus to construct simultaneously two solutions hat have a trend to get closer with time. Each solution being driven by a Brownian motion, a coupling of two solutions then follows from a coupling of the Brownian motions. The nature of the coupling therefore depends on the choice of Brownian motions. The coupling is said to be synchronous if the Brownian motions are the same: it is a coupling for which the Brownian noise cancels out in the infinitesimal evolution of the difference. In that case the contraction of a distance between the processes can only be induced
by the drift component, as in \cite{cattiaux2006probabilisticapproachgranularmedia,BolleyGuillinMalrieu}. Such a contraction holds under the strong convexity conditions.\newline

Nevertheless, in more general cases, the calculation of the evolution of the difference shows that there is still some deterministic contraction on a subspace: see, among others, \cite{eberle2018couplingsquantitativecontractionrates,guillin2021convergenceratesvlasovfokkerplanckequation,kazeykina2023ergodicityunderdampedmeanfieldlangevin}. We can therefore use a synchronous coupling in the vicinity of this subspace. Outside of this subspace, it is necessary to make use of the noise to get the processes closer together, at least in the direction orthogonal to this space. In order to maximize the variance of this noise, we then use a so-called reflection coupling, which consists of taking opposing Brownian motions in the direction of space given by the difference process, and synchronous in the orthogonal direction.\newline

The granular media and the Vlasov-Fokker-Planck dynamics have attracted a lot of attention, usually named McKean-Vlasov diffusions, these past twenty-five years, by means of a \textit{stochastic interpretation} and \textit{synchronous couplings} as in \cite{cattiaux2006probabilisticapproachgranularmedia,BolleyGuillinMalrieu}  or the recent by \textit{reflection couplings}
as in \cite{guillin2021convergenceratesvlasovfokkerplanckequation} enabling to get rid of \textit{convexity conditions}, but limited to \textit{small interactions}. Remark however that \textit{small interactions} are natural to get \textit{uniform in time propagation of chaos} as for large interactions the non linear limit equation may have several stationary measures (see \cite{tugaut2009nonuniquenessstationarymeasures} for example).\newline

Concerning the mean-field Langevin dynamics, propagation of chaos in finite time for the first and second orders McKean-Vlasov dynamics given by \cref{MLD} and \cref{KineticMLD} is relatively easy to show using \textit{synchronous coupling}: for example, regarding \textit{polynomial interaction}, see \cite[Theorem 5.1]{assadeck2023exponentialergodicitymckeanvlasovsde}. The bound obtained by this method, however, generally tends to infinity when the time interval extends to infinity. Besides, the dynamics may possess \textit{multiple invariant measures}, so uniform in time convergence can not be expected without additional assumptions, or a more general definition of convergence itself (e.g. \textit{convergence modulo symmetries}). By assuming the \textit{functional convexity of the energy}, uniform in time propagation of chaos for the mean-field Langevin dynamics and the kinetic mean-field Langevin dynamics have been proven in \cite{chen2023uniformintimepropagationchaosmean,Chen_2024}. In \cite{kazeykina2023ergodicityunderdampedmeanfieldlangevin}, the authors consider a \textit{flat energy} which could
be \textit{non-convex} but is essentially with \textit{small (nonlinear) dependence}, and aim at obtaining an exponential contraction result. The authors mainly borrow the tools developed in  \cite{eberle2017quantitativeharristypetheorems,eberle2018couplingsquantitativecontractionrates}, where they initiate the \textit{reflection-synchronous coupling technique}, further validate it in the study of the Langevin dynamics with a general non-convex potential, and make
the point that the technique offers significant flexibility for additional development. But note
that \cite{eberle2018couplingsquantitativecontractionrates} is not concerned with \textit{mean-field interaction} and the \textit{rate} found there is \textit{dimension
dependent}. In \cite{kazeykina2023ergodicityunderdampedmeanfieldlangevin}, the authors design a new metric involving a quadratic form (\cite[Section 4.4.2]{kazeykina2023ergodicityunderdampedmeanfieldlangevin} to obtain the contraction when the coupled particles are far away, and as a result obtain a \textit{dimension-free convergence rate}. The construction of the quadratic form shares some flavor
with the argument in \cite{BolleyGuillinMalrieu}. Notably, their construction helps to capture
the \textit{optimal rate} in the area of interest (\cite[Remark 4.16]{kazeykina2023ergodicityunderdampedmeanfieldlangevin}), so may be more intrinsic. Notice that most of the articles concerning the ergodicity of underdamped Langevin dynamics obtain
the convergence rates depending on the dimension, and in particular very few allow both non-convex potential and the mean-field interaction. In \cite{BolleyGuillinMalrieu}, the result holds in the \textit{noncompact case} with \textit{small Lipschitz interaction}, and can be seen as a first attempt to deal with more general case. Moreover it shows \textit{existence} and \textit{uniqueness} of the \textit{equilibrium measure}, and in particular does not use its \textit{explicit expression} (which is unknown in our
broader situation). It is also not only a result on the convergence to equilibrium, but also a \textit{stability
result} of all solutions. Let us finally note that it is based on the \textit{natural stochastic interpretation} and a simple coupling argument, and does not need any \textit{hypoelliptic regularity} property of the solutions.\newline

Let us outline the main objectives of the paper. First, sufficient conditions (involving \textit{small interaction}) are provided regarding the macroscopic Hamiltonian $H$ and its \textit{flat derivatives} to ensure the \textit{contraction of the invariant measure map} $\Phi$ for the overdamped mean-field Langevin process (see \cref{mainthm}): this is done using \textit{Lipschitz properties} on the \textit{Poisson equation} and the \textit{Kantorovitch-Rubinstein duality relation}. Secondly, using \textit{synchronous coupling}, we demonstrate \textit{exponential contraction in the Wasserstein} $\mathcal{W}_2$ metric for both overdamped and underdamped mean-field Langevin dynamics (see \cref{W2contraction} and \cref{kineticw2contraction}). Next, we deduce \textit{uniform-in-time $L^2$-propagation of chaos} with a \textit{Fournier–Guillin bound} (see \cref{uniformpropachaosmeanfieldLangevin} and \cref{uniformpropachaoskineticmeanfieldLangevin}). In the kinetic case, general drift conditions are allowed
in the velocity component resulting in a \textit{non-explicit representation of a unique invariant measure}. The results are applied to several concrete and original examples (see \cref{sectionmeanfieldLangevinExamples}).\newline

\textbf{\underline{Main contributions}} In the continuity of \cite{assadeck2023exponentialergodicitymckeanvlasovsde,BolleyGuillinMalrieu}, we are interested in problems of contraction and uniform in time propagation of chaos by synchronous coupling for mean-field Langevin dynamics (over- and under-damped).

\begin{enumerate}
\item\texttt{Over-damped dynamics} 
\begin{itemize}
\item First, as a generalization of the proof of \cite[Proposition 1.1]{assadeck2023exponentialergodicitymckeanvlasovsde}, we show that the invariant measure map is a strict contraction of a metric space (\cref{mainthm}): existence of at most one equilibrium on this space. 
\item Second, by synchronous coupling, we show that there exists a unique equilibrium and that we have exponential convergence towards it (\cref{W2contraction}).
\item Finally, we show a result of propagation of chaos $L^2$ with a bound which can be chosen uniform in time (\cref{uniformpropachaosmeanfieldLangevin}).
\end{itemize}
\item\texttt{under-damped dynamics}
\begin{itemize}
\item First, as a generalization of \cite[Theorem 1]{BolleyGuillinMalrieu}, we show that there exists a unique kinetic equilibrium and that we have exponential convergence towards it (\cref{kineticw2contraction}).
\item Finally, as a generalization of \cite[Theorem 2]{BolleyGuillinMalrieu}, we show a result of uniform in time propagation of chaos $L^2$ (\cref{uniformpropachaoskineticmeanfieldLangevin}).
\end{itemize}
\end{enumerate}

\textbf{\underline{Plan of the paper}}  Let us finish this introduction by the plan of the paper. In \cref{sectionmeanfieldSettings}, we give the settings and assumptions of the paper. The main theorems are stated in \cref{mainresultsfirstorder,mainresultskinetic}. The comments and examples are given in \cref{comments,sectionmeanfieldLangevinExamples}. We discuss the connection with partial differential equations in \cref{sectionmeanfieldLangevinconnectionPDEs}. \cref{sectionproofMainTheoremsMeanfieldLangevin} focuses on the proofs of the main results of the paper. And we end the paper with the bibliographical references.

\section{Assumptions and main theorems}\label{sectionmeanfieldLangevinMainTheorem}
This section focuses on the main hypotheses, the main results and the examples.

\subsection{Settings and assumptions}\label{sectionmeanfieldSettings}
We give an informal preview of our settings and assumptions.\newline

\textbf{\underline{Settings}}
Let $d\geq1$ be an integer. In this paper, we consider mean field energy functionals $H$ that are defined on the space of probability measures on $\mathbb{R}^{d}$. More precisely,
such functionals are mappings ($p\geq0$)
\begin{equation}
    H:\quad\mathcal{P}_{p}(\mathbb{R}^{d})\longrightarrow\mathbb{R}
\end{equation}
where $\mathcal{P}_{p}(\mathbb{R}^{d})$ denotes the set of probability measures of finite $p$-moment equipped with its weak topology. We denote $\mathcal{W}_{p}$ the Wasserstein $p$-metric defined by
\begin{equation}
   \forall (\mu,\nu)\in\mathcal{P}_{p}(\mathbb{R}^{d})\times \mathcal{P}_{p}(\mathbb{R}^{d}),\quad\mathcal{W}_{p}(\mu,\nu):=\inf_{X\sim\mu,Y\sim\nu}(\mathbb{E}[\|X-Y\|^{p}])^{\min(\frac{1}{p},1)}.
\end{equation}
We denote by $\mathcal{D}_{m}H$ and $\mathcal{D}^2_{m}H$ the first and second orders intrinsic derivatives of $H$: the first order intrinsic derivative is a mapping from $\mathcal{P}_{p}(\mathbb{R}^{d})\times\mathbb{R}^{d}$ into $\mathbb{R}^{d}$; and the second order intrinsic derivative is a mapping from $\mathcal{P}_{p}(\mathbb{R}^{d})\times\mathbb{R}^{d}\times\mathbb{R}^{d}$ into $\mathbb{R}^{d\times d}$ and will be defined by
\begin{align}\label{intrinsecderiv}
\mathcal{D}_{m}H(\mu,x)&:=\nabla_{x}\frac{\delta H}{\delta m}(\mu,x);\\
[\mathcal{D}_{m}H]_{1,\infty}&:=\sup_{\mu}[\mathcal{D}_{m}H(\mu,\cdot)]_{1}\textnormal{ (Uniform Lipschitz constant)};\\
\mathcal{D}^2_{m}H(\mu,x,y)&:=\nabla^{2}_{x,y}\frac{\delta^2 H}{\delta m^2}(\mu,x,y);\\
\|\mathcal{D}^{2}_{m}H\|_{\mathbf{op},\infty}&:=\sup_{\mu,x,y}\|\mathcal{D}^2_{m}H(\mu,x,y)\|_{2};
\end{align}
where $\frac{\delta H}{\delta m}$ is called (first) \textit{flat derivative} of $H$ and is a continuous mapping from $\mathcal{P}_{p}(\mathbb{R}^{d})\times\mathbb{R}^{d}$ into $\mathbb{R}$ satisfying the two following properties :
\begin{itemize}
\item For all $K\subset\mathcal{P}_{p}(\mathbb{R}^{d})$ compact,
\begin{equation}
    \exists C_{K}\geq0\quad\forall x\in\mathbb{R}^{d},\quad\sup_{\mu\in K}\bigg|\frac{\delta H}{\delta m}(\mu,x)\bigg|\leq C_{K}(1+\|x\|^p);
\end{equation}
\item\begin{align}\label{deffirstflatderivative}
    H(\mu_{1})-H(\mu_{0})&=\int_{0}^{1}\int_{\mathbb{R}^d}\frac{\delta H}{\delta m}(t\mu_{1}+(1-t)\mu_{0},x)(\mu_{1}-\mu_{0})(dx)dt\\
    \Longleftrightarrow\frac{d}{dt}H(t\mu_{1}+(1-t)\mu_{0})&=\int_{\mathbb{R}^d}\frac{\delta H}{\delta m}(t\mu_{1}+(1-t)\mu_{0},x)(\mu_{1}-\mu_{0})(dx).
\end{align}
\end{itemize}

\textbf{\underline{Assumptions}}
We impose the following assumptions on the energy $H$.
\begin{assumption1}[\texttt{Over-damped dynamics and $\Phi$-contraction}]\label{assumptionPhicontraction}
    The functional $H$ admits first and second order flat derivatives
    \begin{align}
        \frac{\delta H}{\delta m}:\quad&\mathcal{P}_{p}(\mathbb{R}^{d})\times\mathbb{R}^{d}\longrightarrow\mathbb{R};\\
        \frac{\delta^2 H}{\delta m^2}:\quad&\mathcal{P}_{p}(\mathbb{R}^{d})\times\mathbb{R}^{d}\times\mathbb{R}^{d}\longrightarrow\mathbb{R};
    \end{align}
    that are jointly continuous and are $\mathcal{C}^{2}$ in the spatial variables.
    we have the following:
    \begin{enumerate}
        \item $(\mu,x,y)\mapsto\|\mathcal{D}^2_{m}H(\mu,x,y)\|_{2}$ is bounded ($\|\mathcal{D}^{2}_{m}H\|_{\mathbf{op},\infty}<+\infty$);
        \item $(\mu,x)\mapsto\lambda_{\nabla\mathcal{D}_{m}H(\mu,x)}$ (smallest eigenvalue) is strictly bounded from below by $\|\mathcal{D}^{2}_{m}H\|_{\mathbf{op},\infty}$.
    \end{enumerate}
\end{assumption1}

\begin{assumption1}[\texttt{Over-damped dynamics and $\mathcal{W}_{2}$-contraction}]\label{assumptionW2contractioncoupling}
    Assume that the functional $H$ admits first and second order flat derivatives that are jointly continuous and are $\mathcal{C}^{2}$ in the spatial variables. Suppose that
    \begin{equation}
        \sup_{\mu\in\mathcal{P}(\mathbb{R}^{d})}\bigg\|\nabla\frac{\delta H}{\delta m}(\mu,0)\bigg\|_{2}<+\infty
    \end{equation}
    and there exists $\lambda>0$ such that
    \begin{enumerate}
        \item for all $(\mu,x,y)\in\mathcal{P}(\mathbb{R}^{d})\times\mathbb{R}^{d}\times\mathbb{R}^{d}$, 
        \begin{equation}
            \bigg\langle\nabla\frac{\delta H}{\delta m}(\mu,x)-\nabla\frac{\delta H}{\delta m}(\mu,y),x-y\bigg\rangle\geq\lambda \|x-y\|^2_{2};
        \end{equation}
        \item the map $(\mu,x,y)\mapsto\|\mathcal{D}^2_{m}H(\mu,x,y)\|_{2}$ is bounded ($\|\mathcal{D}^{2}_{m}H\|_{\mathbf{op},\infty}<+\infty$) and $\lambda>\|\mathcal{D}^{2}_{m}H\|_{\mathbf{op},\infty}$.
    \end{enumerate}
\end{assumption1}

\begin{assumption1}[\texttt{Under-damped dynamics and $\mathcal{W}_{2}$-contraction}]\label{assumptionkineticW2contractioncoupling}
    Assume that
\begin{itemize}
    \item $[A]_{1}<+\infty$, $B:p\mapsto\lambda_{B}p+D(p)$ with $\lambda_{B}>[A]_{1}$ and $[D]_{1}<+\infty$.
    \item There exists $\lambda_{A}\geq0$ such that for all $(v,w)\in\mathbb{R}^{d}\times\mathbb{R}^{d}$, 
    \begin{equation}
        \langle A(v)-A(w),v-w\rangle\geq\lambda_{A}|v-w|^2.
    \end{equation}
    \item $[\mathcal{D}_{m}H]_{1,\infty}<+\infty$ and $\|\mathcal{D}^{2}_{m}H\|_{\mathbf{op},\infty}<+\infty$.
\end{itemize}
\end{assumption1}

\subsection{Main results for the first-order case}\label{mainresultsfirstorder}
\begin{theorem1}[\texttt{Contraction of $\Phi$}]\label{mainthm}
    If the energy functional $H$ satisfies \cref{assumptionPhicontraction}, then 
    \begin{equation}
        \forall p\geq0,\quad\Phi:\quad\mathcal{P}_{p}(\mathbb{R}^{d})\longrightarrow\mathcal{P}_{p}(\mathbb{R}^{d}),
    \end{equation}
    and for
    \begin{equation}
    \underline{\lambda}:=\inf_{\mu,x}\lambda_{\nabla\mathcal{D}_{m}H(\mu,x)}\quad(>\|\mathcal{D}^{2}_{m}H\|_{\mathbf{op},\infty}),
\end{equation}
for all $p\geq1$, $\Phi$ defines a strict contraction of $(\mathcal{P}_{p}(\mathbb{R}^{d}),\mathcal{W}_{1})$:
\begin{equation}
    \forall (\mu,\nu)\in\mathcal{P}_{p}(\mathbb{R}^{d})\times\mathcal{P}_{p}(\mathbb{R}^{d}),\quad\mathcal{W}_{1}(\Phi(\mu),\Phi(\nu))\leq\frac{\|\mathcal{D}^{2}_{m}H\|_{\mathbf{op},\infty}}{\underline{\lambda}}\mathcal{W}_{1}(\mu,\nu).
\end{equation}
In this case, $\Phi$ admits at most one fixed point on $\mathcal{P}_{p}(\mathbb{R}^{d})$: if $p=1$, there is a unique fixed point.
\end{theorem1}
\begin{proof}
    See \cref{thm1}.
\end{proof}

\begin{theorem1}[\texttt{$\mathcal{W}_{2}$-Exponential contraction by coupling}]\label{W2contraction}
If the energy functional $H$ satisfies \cref{assumptionW2contractioncoupling}, for all $(X_{t})_{t\geq0}$ and $(Y_{t})_{t\geq0}$ solutions of the mean field Langevin equation \cref{MLD} with initial conditions $\mu^{X }_{0}\in\mathcal{P}_{2}(\mathbb{R}^{d})$ and $\mu^{Y}_{0}\in\mathcal{P}_{2}(\mathbb{R}^{d})$, we have
\begin{itemize}
    \item for all $t\geq0$, $\mu^{X}_{t}\in\mathcal{P}_{2}(\mathbb{R}^{d})$;
    \item for all $t\geq0$,
    \begin{equation}
        \mathcal{W}_{2}(\mu^{X}_{t},\mu^{Y}_{t})\leq e^{(\|\mathcal{D}^{2}_{m}H\|_{\mathbf{op},\infty}-\lambda)t}\mathcal{W}_{2}(\mu^{X}_{0},\mu^{Y}_{0}).
    \end{equation}
\end{itemize}
In particular, there exists a unique invariant measure $\mu_{\infty}\in\mathcal{P}_{2}(\mathbb{R}^{d})$ and
\begin{equation}
    \forall t\geq0,\quad\mathcal{W}_{2}(\mu^{X}_{t},\mu_{\infty})\leq e^{(\|\mathcal{D}^{2}_{m}H\|_{\mathbf{op},\infty}-\lambda)t}\mathcal{W}_{2}(\mu^{X}_{0},\mu_{\infty}).
\end{equation}
\end{theorem1}
\begin{proof}
    See \cref{proofw2contraction}.
\end{proof}

\begin{theorem1}[\texttt{Uniform in time propagation of chaos for mean field Langevin dynamics}]\label{uniformpropachaosmeanfieldLangevin}
    If the energy functional $H$ satisfies \cref{assumptionW2contractioncoupling}, for $(\widetilde{X}^{(n)})_{n\geq1}$ a sequence of independent copies of the mean field Langevin process $(X_{t})_{t\geq0}$ synchronously coupled to particles of the over-damped mean field Langevin system, we have
    \begin{align}
        \forall T\geq0\quad\forall t\in [0,T],\quad\mathbb{E}[\|X^{(1)}_{t}-\widetilde{X}^{(1)}_{t}\|^2]&\leq\begin{cases}
            \frac{\alpha_{1}(T)}{\beta_{1}}(e^{\beta_{1} t}-1)\quad&\textnormal{if}\quad\beta_{1}\neq0;\\
            \alpha_{1}(T) t\quad&\textnormal{if}\quad\beta_{1}=0.
        \end{cases}\\
        \alpha_{1}(T)&:=\|\mathcal{D}^2_{m}H\|_{\mathbf{op},\infty}C(d)\bigg(\sup_{t\in[0,T]}\bigg\langle\mu^{X}_{t},\|\cdot\|^2_{2}\bigg\rangle\bigg)\delta_{d}(n);\\
        \delta_{d}(n)&:=\begin{cases}
             \frac{1}{\sqrt{n}}\quad&\textnormal{if}\quad d<4;\\
             \frac{\ln(n+1)}{\sqrt{n}}\quad&\textnormal{if}\quad d=4;\\
             \frac{1}{n^{\frac{2}{d}}}\quad&\textnormal{if}\quad d>4;
         \end{cases}\\
        \beta_{1}&:=3\|\mathcal{D}^2_{m}H\|_{\mathbf{op},\infty}-2\lambda;\\
        \alpha&:=2d+2\sup_{\mu}\bigg\|\nabla\frac{\delta H}{\delta m}(\mu,0)\bigg\|_{2};\\
        \beta&:=-2\bigg(\lambda-\sup_{\mu}\bigg\|\nabla\frac{\delta H}{\delta m}(\mu,0)\bigg\|_{2}\bigg).
    \end{align}
    Moreover, this upper bound is uniform in time if, and only if, $\beta<0$ and $\beta_{1}<0$ :
    \begin{equation}
        \sup_{t\geq0}\mathbb{E}[\|X^{(1)}_{t}-\widetilde{X}^{(1)}_{t}\|^2]\leq\frac{\|\mathcal{D}^2_{m}H\|_{\mathbf{op},\infty}C(d)}{2\lambda-3\|\mathcal{D}^2_{m}H\|_{\mathbf{op},\infty}}\bigg(\bigg\langle\mu^{X}_{0},\|\cdot\|^2_{2}\bigg\rangle-\frac{\alpha}{\beta}\bigg)\delta_{d}(n).
    \end{equation}
\end{theorem1}

\begin{proof}
    See \cref{proofuniformpropachaosmeanfieldLangevin}.
\end{proof}

\subsection{Main results for the kinetic case}\label{mainresultskinetic}
\begin{theorem1}[\texttt{Kinetic contraction by coupling}]\label{kineticw2contraction}
If the fields $A$ and $B$ and the energy functional $H$ satisfy \cref{assumptionkineticW2contractioncoupling}, for all positive $[A]_{1}$, $\lambda_{A}$ and $\lambda_{B}$, there exists a positive constant $c$ such that, if  
    \begin{equation}
        0\leq\min([D]_{1},[\mathcal{D}_{m}H]_{1,\infty},\|\mathcal{D}^{2}_{m}H\|_{\mathbf{op},\infty})\leq\max([D]_{1},[\mathcal{D}_{m}H]_{1,\infty},\|\mathcal{D}^{2}_{m}H\|_{\mathbf{op},\infty})<c,
    \end{equation}
    then there exist positive constants $C$ and $C'$ such that, 
    for all $(P_{t},V_{t})_{t\geq0}$ and $(\overline{P}_{t},\overline{V}_{t})_{t\geq0}$ two solutions of the kinetic mean-field Langevin equation \cref{KineticMLD} with initial conditions $\mathbb{P}_{(P_{0},V_{0})}\in\mathcal{P}_{2}(\mathbb{R}^{d}\times\mathbb{R}^{d})$ and $\mathbb{P}_{(\overline{P}_{0},\overline{V}_{0})}\in\mathcal{P}_{2}(\mathbb{R}^{d}\times\mathbb{R}^{d})$, we have
    \begin{equation}
        \forall t\geq0,\quad\mathcal{W}_{2}(\mathbb{P}_{(P_{t},V_{t})},\mathbb{P}_{(\overline{P}_{t},\overline{V}_{t})})\leq C'e^{-Ct}\mathcal{W}_{2}(\mathbb{P}_{(P_{0},V_{0})},\mathbb{P}_{(\overline{P}_{0},\overline{V}_{0})}).
    \end{equation}
    Moreover, there exists a unique stationary solution $\xi_{\infty}$ and
    \begin{equation}
        \forall t\geq0,\quad\mathcal{W}_{2}(\mathbb{P}_{(P_{t},V_{t})},\xi_{\infty})\leq C'e^{-Ct}\mathcal{W}_{2}(\mathbb{P}_{(P_{0},V_{0})},\xi_{\infty}).
    \end{equation}
\end{theorem1}

\begin{proof}
    See \cref{proofkineticw2contraction}.
\end{proof}

\begin{theorem1}[\texttt{Uniform in time propagation of chaos for kinetic mean field Langevin dynamics}]\label{uniformpropachaoskineticmeanfieldLangevin}
    Let $((\overline{P}^{(n)},\overline{V}^{(n)}))_{n\geq 1}$ a sequence of independent copies of the kinetic mean field Langevin process $((P_{t},V_{t}))_{t\geq0}$ synchronously coupled to particles of the under-damped mean field Langevin system with initial datum $(P^{(i)}_{0},V^{(i)}_{0})$ for $i\in\{1,\ldots,n\}$ be $n$ independent $\mathbb{R}^{2d}$-valued random variables with law $f_{0}\in\mathcal{P}_{2}(\mathbb{R}^{2d})$. If the fields $A$ and $B$ and the energy functional $H$ satisfy \cref{assumptionkineticW2contractioncoupling}, for all positive $[A]_{1}$, $\lambda_{A}$ and $\lambda_{B}$, there exists a positive constant $c$ such that, if  
    \begin{equation}
        0\leq\min([D]_{1},[\mathcal{D}_{m}H]_{1,\infty},\|\mathcal{D}^{2}_{m}H\|_{\mathbf{op},\infty})\leq\max([D]_{1},[\mathcal{D}_{m}H]_{1,\infty},\|\mathcal{D}^{2}_{m}H\|_{\mathbf{op},\infty})<c,
    \end{equation}
    then there exists a positive constant $C$, independent of $n$ , such that for $i\in\{1,\ldots,n\}$
    \begin{equation}
        \sup_{t\geq0}\mathbb{E}\bigg[|P^{(i)}_{t}-\overline{P}^{(i)}_{t}|^{2}+|V^{(i)}_{t}-\overline{V}^{(i)}_{t}|^{2}\bigg]\leq C\begin{cases}
             \frac{1}{\sqrt{n}}\quad&\textnormal{if}\quad d<4;\\
             \frac{\ln(n+1)}{\sqrt{n}}\quad&\textnormal{if}\quad d=4;\\
             \frac{1}{n^{\frac{2}{d}}}\quad&\textnormal{if}\quad d>4.
         \end{cases}
    \end{equation}
    Here the constant $C$ depends only on the coefficients of the equation and the second moment of $f_{0}$. 
\end{theorem1}

\begin{proof}
    See \cref{proofuniformpropachaoskineticmeanfieldLangevin}.
\end{proof}

\subsection{Comments on drift and assumptions}\label{comments}
In \cref{W2contraction}, \cref{uniformpropachaosmeanfieldLangevin}, \cref{kineticw2contraction} and \cref{uniformpropachaoskineticmeanfieldLangevin}, we can replace the intrinsic derivative by a more general drift
\begin{equation}
    b:\quad\mathbb{R}^{d}\times\mathcal{P}_{p}(\mathbb{R}^{d})\longrightarrow\mathbb{R}^{d}.
\end{equation}
In the first order case, we can replace \cref{assumptionW2contractioncoupling} by the following conditions on the new drift $b$
\begin{align}
    \sup_{\mu}b(0,\mu)&<+\infty;\\
    \exists\lambda>\beta\geq0\quad\forall (x,y,\mu,\nu),\quad\langle b(x,\mu)-b(y,\nu),x-y\rangle&\geq\lambda|x-y|^{2}-\beta|x-y|\mathcal{W}_{p}(\mu,\nu).
\end{align}
As for the proofs, it is enough to mimic those with the intrinsic derivative (\cref{proofw2contraction,proofuniformpropachaosmeanfieldLangevin}) which is why we do them without loss of generality with the intrinsic derivative.

In kinetic case, we can replace \cref{assumptionkineticW2contractioncoupling} by the same conditions on the velocity field $A$ and the external confinement field $B$ and the following Lipschtz conditions on the new drift $b$
\begin{align}
    \exists\gamma_{1},\gamma_{2}\geq0,\quad\quad\forall (x,y,\mu,\nu),\quad|b(x,\mu)-b(y,\mu)|&\leq\gamma_{1}|x-y|;\\
    |b(x,\mu)-b(x,\nu)|&\leq\gamma_{2}\mathcal{W}_{p}(\mu,\nu).
\end{align}
As in the previous case, we carry out the proofs with the intrinsic derivative without loss of generality.

Recall that results regarding exponential contraction have been established via synchronous coupling with the Wasserstein metric for McKean-Vlasov SDEs featuring constant diffusion coefficients and general drift terms. It would be of interest to derive these results using alternative coupling arguments and/or different metrics: in the case of mean-field kinetic Langevin dynamics, this was done in \cite{kazeykina2023ergodicityunderdampedmeanfieldlangevin}. Furthermore, analogous to \cite[Theorem 1.5]{BolleyGuillinMalrieu}, we can demonstrate a concentration inequality. Finally, given that the results on uniform-in-time propagation of chaos rely on a Fournier-Guillin bound of order $O(1/\sqrt{n})$, it would be worthwhile to revisit the proofs to improve this bound to $O(1/n)$, similar to the result in \cite[Theorem 1.2]{BolleyGuillinMalrieu}.

\subsection{Examples}\label{sectionmeanfieldLangevinExamples}
We give concrete and original examples of flat energies for which our results are satisfied.\newline

\textbf{\underline{Two-body interaction: granular media flat energy}}
In this case, we consider the following flat energy
\begin{equation}
    H_{VW}:\quad\mu\longmapsto\int V(x)\mu(dx)+\frac{1}{2}\int W(x-y)\mu(dx)\mu(dy)
\end{equation}
with $V\in\mathcal{C}^{2}(\mathbb{R}^{d},\mathbb{R})$ called confinement potential and $W\in\mathcal{C}^{2}(\mathbb{R}^{d},\mathbb{R})$ called interaction potential such that for all $x\in\mathbb{R}$, $W(-x)=W(x)$. 

\begin{proposition1}[\texttt{Two-body interaction and \cref{assumptionPhicontraction}}]\label{proposition1Example1}
    If $\|\nabla^2W\|_{\mathbf{op},\infty}<+\infty$ and
    \begin{equation}
        \inf_{(\mu,x)}\lambda_{\nabla^2V(x)+\nabla^{2}W*\mu(x)}>\|\nabla^2W\|_{\mathbf{op},\infty},
    \end{equation}
    $H_{VW}$ satisfies \cref{assumptionPhicontraction}.
\end{proposition1}
\begin{proof}
    See \cref{proofpropo1example1}.
\end{proof}
\begin{remark}[\texttt{About one-dimensional example}]
    In this case, the contraction conditions become
\begin{align*}
    \|W''\|_{\infty}&<+\infty;\\
    \forall (\mu,x),\quad V''(x)+\int W''(x-y)\mu(dy)&>\|W''\|_{\infty}.
\end{align*}
We have
\begin{align*}
    V''(x)+\int W''(x-y)\mu(dy)\geq V''(x)-\|W''\|_{\infty}.
\end{align*}
In particular, if for all $x\in\mathbb{R}$, $V''(x)>2\|W''\|_{\infty}$, we deduce that
\begin{align*}
    \forall (\mu,x),\quad V''(x)+\int W''(x-y)\mu(dy)>\|W''\|_{\infty}.
\end{align*}
For example, we can choose the following polynomial confinement potentials:\\

$\rhd$ $V(x):=ax^{2}+bx+c$ with $a,b$ and $c$ constants such that $a>\|W''\|_{\infty}$ and $(b,c)\in\mathbb{R}^{2}$.\\

$\rhd$ $V(x):=ax^{4}+bx^{3}+cx^{2}+dx+e$ with $a,b,c,d$ and $e$ constants such that $a>0$, $c-\frac{3b^{2}}{8a}>\|W''\|_{\infty}$ and $(d,e)\in\mathbb{R}^{2}$.
\end{remark}
\begin{proposition1}[\texttt{Two-body interaction and \cref{assumptionW2contractioncoupling}}]\label{proposition2Example1}
    If $\|\nabla W\|_{2,+\infty}<+\infty$, $\|\nabla^2W\|_{\mathbf{op},\infty}<+\infty$ and there exists $\lambda>0$ such that for all $(\mu,x,y)$,
    \begin{align}
        \bigg\langle\nabla V(x)-\nabla V(y)+\nabla W*\mu(x)-\nabla W*\mu(y),x-y\bigg\rangle&\geq\lambda\|x-y\|^2_{2};\\
        \lambda&>\|\nabla^2W\|_{\mathbf{op},\infty};
    \end{align}
    $H_{VW}$ satisfies \cref{assumptionW2contractioncoupling}.
\end{proposition1}

\begin{proof}
    See \cref{proofpropo2example1}.
\end{proof}

\begin{remark}
    If for all (x,y),
    \begin{align}
        \bigg\langle\nabla V(x)-\nabla V(y),x-y\bigg\rangle&\geq\lambda_{V}\|x-y\|^2_{2};\\
        \bigg\langle\nabla W(x)-\nabla W(y),x-y\bigg\rangle&\geq0;
    \end{align}
    we have
    \begin{align}
        \bigg\langle\nabla V(x)-\nabla V(y)+\nabla W*\mu(x)-\nabla W*\mu(y),x-y\bigg\rangle&=\bigg\langle\nabla V(x)-\nabla V(y),x-y\bigg\rangle\\
        &+\int\bigg\langle\nabla W(x-z)-\nabla W(y-z),x-z-(y-z)\bigg\rangle\mu(dz)\nonumber\\
        &\geq\lambda_{V}\|x-y\|^2_{2}.\nonumber
    \end{align}
    Note that, here, the obtained conditions match with conditions considered in \cite{malrieu:hal-01282602}.
\end{remark}

\textbf{\underline{Many-body interaction: polynomial flat energy}}
In this case, we consider the following flat energy
\begin{equation}
    H^{N}_{VW}:\quad\mu\longmapsto\int Vd\mu+\sum_{k=2}^{N}\int W^{(k)}d\mu^{\otimes k},
\end{equation}
where $\forall k\in\{2,\ldots,N\}$, $W^{(k)}$ is a \textit{symmetric interaction potential} between $k$ particles and $N$ represents the number of such potentials. 

\begin{proposition1}[\texttt{Many-body interaction and \cref{assumptionPhicontraction}}]\label{proposition1Example2}
    If for all $k\in\{2,\ldots,N\}$, $\|\nabla^2_{x_{1},x_{2}}W^{(k)}\|_{\mathbf{op},\infty}<+\infty$ and
    \begin{equation}
        \inf_{(\mu,x)}\lambda_{\nabla^2V(x)+\sum_{k=2}^{N}k\nabla^2_{x_{1}}W^{(k)}*\mu^{\otimes k-1}(x)}>\sum_{k=2}^{N}k(k-1)\|\nabla^2_{x_{1},x_{2}}W^{(k)}\|_{\mathbf{op},\infty},
    \end{equation}
    $H^{N}_{VW}$ satisfies \cref{assumptionPhicontraction}.
\end{proposition1}
\begin{proof}
    See \cref{proofpropo1example2}.
\end{proof}

\begin{proposition1}[\texttt{Many-body interaction and \cref{assumptionW2contractioncoupling}}]\label{proposition2Example2}
    If for all $k\in\{2,\ldots,N\}$, $\|\nabla_{x_{1}} W^{(k)}(0,\cdot)\|_{2,+\infty}<+\infty$, $\|\nabla^2_{x_{1},x_{2}}W^{(k)}\|_{\mathbf{op},\infty}<+\infty$ and there exists $\lambda>0$ such that for all $(\mu,x,y)$,
    \begin{align}
        &\bigg\langle\nabla V(x)-\nabla V(y)+\sum_{k=2}^{N}k\bigg(\nabla_{x_{1}} W^{(k)}*\mu^{\otimes k-1}(x)-\nabla_{x_{1}} W^{(k)}*\mu^{\otimes k-1}(y)\bigg),x-y\bigg\rangle\geq\lambda\|x-y\|^2_{2};\\
        &\lambda>\sum_{k=2}^{N}k(k-1)\|\nabla^2_{x_{1},x_{2}}W^{(k)}\|_{\mathbf{op},\infty};
    \end{align}
    $H^{N}_{VW}$ satisfies \cref{assumptionW2contractioncoupling}.
\end{proposition1}

\begin{proof}
    See \cref{proofpropo2example2}.
\end{proof}

\begin{remark} For further details on many-body interactions and U-statistics, see \cite{aboubacrineassadek:tel-05460067,assadeck2023exponentialergodicitymckeanvlasovsde}. Note that the \textit{granular-media flat energy} is a particular case of the polynomial energy:
\begin{align}
    H_{VW}(\mu)&=\int V(x)\mu(dx)+\int W^{(2)}(x,y)\mu(dx)\mu(dy);\\
    W^{(2)}(x,y)&:=\frac{1}{2}W(x-y).
\end{align}
We deduce from this that \cref{proposition1Example2} and \cref{proposition2Example2} generalize \cref{proposition1Example1} and \cref{proposition2Example1}.
\end{remark}

\textbf{\underline{Internal flat energy}}
In this case, we consider the following flat energy
\begin{equation}
    H_{\psi}:\quad\mu\longmapsto\psi(\langle\mu,W\rangle)
\end{equation}
with $\psi\in\mathcal{C}^{2}(\mathbb{R},\mathbb{R})$ and $W\in\mathcal{C}^{2}(\mathbb{R}^{d},\mathbb{R})$. 
\begin{proposition1}[\texttt{Internal flat energy and \cref{assumptionPhicontraction}}]\label{proposition1Example3}
    If $\sup_{\mu}|\psi''(\langle\mu,W\rangle)|<+\infty$, $\|\nabla^2 (W\otimes W)\|_{\mathbf{op},\infty}<+\infty$ and
    \begin{equation}
        \inf_{(\mu,x)}\lambda_{\psi'(\langle\mu,W\rangle)\nabla^{2}W(x)}>\|\nabla^2(W\otimes W)\|_{\mathbf{op},\infty}\sup_{\mu}|\psi''(\langle\mu,W\rangle)|,
    \end{equation}
    $H_{\psi}$ satisfies \cref{assumptionPhicontraction}.
\end{proposition1}
\begin{proof}
    See \cref{proofpropo1example3}.
\end{proof}

\begin{proposition1}[\texttt{Internal flat energy and \cref{assumptionW2contractioncoupling}}]\label{proposition2Example3}
    If $\sup_{\mu}|\psi'(\langle\mu,W\rangle)|<+\infty$, $\sup_{\mu}|\psi''(\langle\mu,W\rangle)|<+\infty$, $\|\nabla^2 (W\otimes W)\|_{\mathbf{op},\infty}<+\infty$ and there exists $\lambda>0$ such that for all $(\mu,x,y)$,
    \begin{align}
        \psi'(\langle\mu,W\rangle)\bigg\langle\nabla W(x)-\nabla W(y),x-y\bigg\rangle&\geq\lambda\|x-y\|^2_{2};\\
        \lambda&>\|\nabla^2(W\otimes W)\|_{\mathbf{op},\infty}\sup_{\mu}|\psi''(\langle\mu,W\rangle)|;
    \end{align}
    $H_{\psi}$ satisfies \cref{assumptionW2contractioncoupling}.
\end{proposition1}

\begin{proof}
    See \cref{proofpropo2example3}.
\end{proof}

\section{Connection with partial differential equations}\label{sectionmeanfieldLangevinconnectionPDEs}
The mean field Langevin dynamics corresponds to  nonlinear Fokker-Planck equation defined on $[0,+\infty[\times\mathbb{R}^{d}$ for the flow of probability measures $(\mu_{t})_{t\geq0}$:
\begin{equation}\label{meanfieldPDE}
    \frac{d}{dt}\mu_{t}=\Delta\mu_{t}+\nabla\cdot(\mu_{t}\mathcal{D}_{m}H(\mu_{t},\cdot)).
\end{equation}
In the kinetic case, we have
\begin{equation}\label{meanfieldkineticPDE}
    \frac{d}{dt}f_{t}+v\cdot\nabla_{p}f_{t}-\mathcal{D}_{m}H(f_{t},p)\cdot\nabla_{v}f_{t}=\Delta_{v}f_{t}+\nabla_{v}\cdot\bigg((A(v)+B(p))f_{t}\bigg).
\end{equation}
The particle system corresponds to a Liouville or Fokker-Planck equation defined on $[0,+\infty[\times\mathbb{R}^{nd}$ for the flow of probability measures defined by
\begin{align}
    \forall t\geq0,\quad\mu^{(n)}_{t}&:=\mathbb{P}_{\mathbf{X}_{t}};\\
    H_{n}(\mathbf{x})&:=nH(\mu_{\mathbf{x}});\\
    \frac{d}{dt}\mu^{(n)}_{t}&=\Delta\mu^{(n)}_{t}+\nabla\cdot(\mu^{(n)}_{t}\nabla H_{n}).
\end{align}
In the kinetic case, we have the kinetic Fokker-Planck equation :
\begin{equation}
    \frac{d}{dt}f^{(n)}_{t}+\mathbf{v}\cdot\nabla_{\mathbf{p}}f^{(n)}_{t}-\nabla H_{n}(\mathbf{p})\cdot\nabla_{\mathbf{v}}f^{(n)}_{t}=\Delta_{\mathbf{v}}f^{(n)}_{t}+\sum_{i=1}^{n}\nabla_{v_{i}}\cdot\bigg((A(v_{i})+B(p_{i}))f^{(n)}_{t}\bigg).
\end{equation}

\begin{theorem1}[\texttt{Characterization of invariant distributions}]\label{caracterisationIM}
Assume that for all $\mu\in\mathcal{P}(\mathbb{R}^{d})$, $x\in\mathbb{R}^{d}\mapsto\frac{\delta H}{\delta m}(\mu,x)$ is differentiable and 
    \begin{equation}
        Z_{\mu}:=\int_{\mathbb{R}^{d}}e^{-\frac{\delta H}{\delta m}(\mu,x)}dx<+\infty.
    \end{equation}
    The following assertions are equivalent:
\begin{enumerate}
    \item $\mu_{\infty}$ is a Maxwellian of the McKean-Vlasov PDE:
    \begin{equation}
        \mu_{\infty}\in\{\nu,\quad\Delta\nu+\nabla\cdot(\nu\mathcal{D}_{m}H(\nu,\cdot))=0\}
    \end{equation}
    \item $\mu_{\infty}$ is a critical point for the free energy:
    \begin{equation}
        \mu_{\infty}\in\bigg\{\nu,\quad\mathcal{D}_{m}H(\nu,\cdot)+\nabla\ln\frac{d\nu}{dx}=0\bigg\}.
    \end{equation}
    \item $\mu_{\infty}$ is a fixed point of $\Phi$:
    \begin{equation}
        \mu_{\infty}\in\{\nu,\quad\Phi(\nu)=\nu\}.
    \end{equation}
\end{enumerate}    
\end{theorem1}
\begin{proof}
    See \cref{thm2}.
\end{proof}

In the kinetic case, we recall that if the velocity and external confinement fields do not derive from potential gradients, the invariant probabilities are not explicit : in this paper, we will use \cite[Lemma 10]{BolleyGuillinMalrieu} (contraction property of a semigroup on a complete metric space and existence of a unique fixed point).

\section{Proof of main theorems}\label{sectionproofMainTheoremsMeanfieldLangevin}

\subsection{Proof of main results of the first order case}
\begin{proof}[Proof of \cref{mainthm}]\label{thm1}
We first check that $\Phi(\mathcal{P}_p(\R^d) )\subset {\cal P}_p(\R^d)$. In \cref{assumptionPhicontraction}, the uniform lower bound of the smallest eigenvalue of the Hessian of the first flat derivative implies that
\begin{equation}
    \exists c^{H}_{1}>0\quad\exists c^{H}_{2}\in\mathbb{R}\quad\forall x\in\mathbb{R}^{d}\quad\forall\mu\in\mathcal{P}_{p}(\mathbb{R}^{d}),\quad\frac{\delta H}{\delta m}(\mu,x)\geq c^{H}_{1}\|x\|^2_{2}+c^{H}_{2}.
\end{equation}
We deduce that
\begin{align}
    \langle\Phi(\mu),\|\cdot\|^{p}\rangle&=\frac{1}{Z_{\mu}}\int_{\mathbb{R}^{d}}\|x\|^{p}e^{-\frac{\delta H}{\delta m}(\mu,x)}dx\\
    &\leq\frac{e^{-c^{H}_{2}}}{Z_{\mu}}\int_{\mathbb{R}^{d}}\|x\|^{p}e^{-c^{H}_{1}\|x\|^2}dx,\nonumber
\end{align}
and this last integral is finite since
\begin{equation}
    \exists M^{H}>0\quad\forall x\in\mathbb{R}^{d},\quad \|x\|^{p}e^{-c^{H}_{1}\|x\|^2}\leq M^{H}e^{-\frac{c^{H}_{1}}{2}\|x\|^2}.
\end{equation}

    Now, let's prove the contraction of $\Phi$. Let $\nu_{0},\nu_{1}\in {\cal P}_p(\R^{d})$ and set $\nu_{t}=(1-t)\nu_{0}+t\nu_{1}$, $t\in[0,1]$. Let $f:\mathbb{R}^d\rightarrow\mathbb{R}$ be a $1$-Lipschitz smooth function. From the very definition of $\Phi$,
\begin{equation}\label{eq:phinut3}
    \langle\Phi(\nu_{t}),f\rangle=\frac{1}{Z_{\nu_{t}}}\int_{\mathbb{R}^{d}}f(x)e^{-\frac{\delta H}{\delta m}(\nu_{t},x)}dx,
\end{equation}
so that setting 
\begin{equation}
    g_t(x):=-\frac{d}{dt}\frac{\delta H}{\delta m}(\nu_{t},x),
\end{equation}
we get
\begin{align}
\frac{d}{dt}\langle\Phi(\nu_{t}),f\rangle&=-\frac{\partial_t Z_{\nu_t}}{Z_{\nu_t}}\langle\Phi(\nu_{t}),f\rangle+\langle\Phi(\nu_{t}),f g_t \rangle\\
&=\langle\Phi(\nu_{t}),f g_t \rangle-\langle\Phi(\nu_{t}),f\rangle\langle\Phi(\nu_{t}),g_t\rangle\nonumber\\
&={\rm Cov}_{\Phi(\nu_{t})}(f,g_t).\nonumber
\end{align}
As the functional $H$ admits first and second-order flat derivatives, we have
\begin{equation}
    g_t(x):=-\int_{\mathbb{R}^{d}}\frac{\delta^2 H}{\delta m^2}(\nu_{t},x,y)(\nu_{1}-\nu_{0})(dy).
\end{equation}
For a probability $\mu$, let ${\cal L}_{\mu}$ be the operator defined on ${\cal C}^2$-functions by 
\begin{equation}\label{eq:operator3}
 {\cal L}_\mu f=\Delta f-\mathcal{D}_{m}H(\mu,\cdot)\cdot\nabla f.
 \end{equation}
 Denoting by $\varphi_t$ be the solution of the Poisson equation $f-\Phi(\nu_{t})(f)={\cal L}_{\nu_t} \varphi_t$ and using that ${\cal L}_{\nu_t}$ is self-adjoint in $L^2(\Phi(\nu_{t}))$, we get
$${\rm Cov}_{\Phi(\nu_{t})}(f,g_t)=\langle \varphi_t, {\cal L}_{\nu_t} g_t\rangle_{L^2(\Phi(\nu_{t}))}=-\langle \nabla \varphi_t, \nabla g_t\rangle_{L^2(\Phi(\nu_{t}))}.$$
Note that for the second equality, we used the fact that for some ${\cal C}^2$-functions $f$ and $g$,  ${\cal L}_{\nu_t} (f.g)= f {\cal L}_{\nu_t} g+g{\cal L}_{\nu_t} f+2\nabla f. \nabla g$. With the help of Cauchy-Schwarz inequality, this leads to 
\begin{align}\label{eq:pluggt3}
|{\rm Cov}_{\Phi(\nu_{t})}(f,g_t)|&\le \|\nabla \varphi_t\|_{L^2(\Phi(\nu_{t}))}\|\nabla g_t\|_{L^2(\Phi(\nu_{t}))}\\
&\le  \|\nabla \varphi_t\|_{\infty}\|\nabla g_t\|_{\infty}.\nonumber
\end{align}

\noindent\textbf{Control of $\|\nabla \varphi_t\|_{\infty}$.} Let $h:\mathbb{R}^d\rightarrow\mathbb{R}$ be a ${\cal C}^1$ function with bounded derivative. Let $\varphi$ be a unique solution of the Poisson equation $h-\langle \Phi(\mu), h\rangle={\cal L}_\mu \varphi$. It is well-known that $\varphi(x)=\int_0^{+\infty} P_t^{\mu} h(x)-\langle\Phi(\mu),h\rangle dt$ where $(P_t^{\mu})_{t\ge0}$ denotes the semi-group associated with ${\cal L}_\mu$ so that under adequate derivation conditions (which will be satisfied in our setting), 
$$\nabla \varphi(x)=\int_0^{+\infty} \nabla P_t^{\mu} h(x) dt.$$
Indeed, it is sufficient to verify the conditions of application of the Lebesgue dominated convergence theorem for the derivation under the integral sign, which is done below.

We have 
$$ \nabla P_t^{\mu} h (x)=\mathbb{E}[\nabla h(X_t^{\mu,x}) \partial_x X_t^{\mu,x}]$$
where $(X_t^{\mu,x})$ denotes the solution starting from $x$ of the SDE associated with ${\cal L}_\mu$ and $(\partial_x X_t^\mu)$ its first variation process solution to 
$$ dY_t=-\nabla\mathcal{D}_{m}H(\mu,X^{\mu,x}_{t}) Y_t dt$$
with $Y_0=I_d$. Under the assumption 
\begin{equation}
    \underline{\lambda}:=\inf_{\mu,x}\lambda_{\nabla\mathcal{D}_{m}H(\mu,x)}>0,
\end{equation}
one easily deduces from a Gronwall argument that for any $z\in\mathbb{R}^d$,
$$ |Y_t z|^2\le e^{-2\underline{\lambda}t}|z|^2$$
which implies that 
$$|\nabla P_t^{\mu} h(x)|\le [h]_1\|\partial_x X_t^{\mu,x}\|_{\mathbf{op}}\le [h]_1 e^{-\underline{\lambda}t}$$
where $[h]_1$ denotes the Lipschitz constant of $h$.
Then,
$$\|\nabla \varphi\|_\infty\le [h]_1\underline{\lambda}^{-1}.$$
In particular, we have
\begin{equation}
    \|\nabla \varphi_t\|_{\infty}\leq[f]_1\underline{\lambda}^{-1}=\frac{1}{\underline{\lambda}}.
\end{equation}

\noindent\textbf{Control of $\|\nabla g_t\|_{\infty}$.} Let us finally focus on $\nabla g_t$. First, one checks that
\begin{equation}
    \nabla g_t(x):=-\int_{\mathbb{R}^{d}}\nabla_{x}\frac{\delta^2 H}{\delta m^2}(\nu_{t},x,y)(\nu_{1}-\nu_{0})(dy).
\end{equation}
Now, denoting by $(X_{\nu_0},X_{\nu_1})$ an optimal coupling of $\nu_0$ and $\nu_1$ for the $1$-Wasserstein distance, one obtains
\begin{align}
\bigg|\nabla_x \left(\int_{\mathbb{R}^{d}}\frac{\delta^2 H}{\delta m^2}(\nu_{t},x,y)(\nu_{1}-\nu_{0})(dy)\right)\bigg|&\le  |\mathbb{E}[\nabla_{x}\frac{\delta^2 H}{\delta m^2}(\nu_{t},x,X_{\nu_{1}}) -\nabla_{x}\frac{\delta^2 H}{\delta m^2}(\nu_{t},x,X_{\nu_{0}})]|\\
&\le \|\mathcal{D}^{2}_{m}H\|_{\mathbf{op},\infty}\mathcal{W}_{1}(\nu_{0},\nu_{1}).\nonumber
\end{align}
Hence, 
$$ \|\nabla g_t\|_{\infty}\le \|\mathcal{D}^{2}_{m}H\|_{\mathbf{op},\infty}\mathcal{W}_{1}(\nu_{0},\nu_{1}),$$
and by \eqref{eq:pluggt3}, we get for any smooth $1$-Lipschitz function
$$|\langle\Phi(\nu_{1}),f\rangle-\langle\Phi(\nu_{0}),f\rangle|\leq\int_{0}^{1}|\frac{d}{dt}\langle\Phi(\nu_{t}),f\rangle| dt\le \frac{\|\mathcal{D}^{2}_{m}H\|_{\mathbf{op},\infty}}{\underline{\lambda}}\mathcal{W}_{1}(\nu_{0},\nu_{1}).$$
Since
\begin{equation}
    \frac{\|\mathcal{D}^{2}_{m}H\|_{\mathbf{op},\infty}}{\underline{\lambda}}<1,
\end{equation}
by a density argument and the Kantorovitch-Rubinstein duality relation, it follows that $\Phi$ is a contraction on $({\cal P}_p(\R^{d}),\mathcal{W}_{1})$.\newline

Concerning the existence of at most one fixed point, we show it by the absurd: assuming that we have more than one, we would contradict the strict contraction since we would have
\begin{equation}
    \frac{\|\mathcal{D}^{2}_{m}H\|_{\mathbf{op},\infty}}{\underline{\lambda}}\geq1.
\end{equation}
For $p=1$, existence follows from Picard's fixed point theorem since $\Phi$ is a strict contraction of $(\mathcal{P}_{1}(\mathbb{R}^{d}),\mathcal{W}_{1})$ which is a complete metric space.
\end{proof}

\begin{proof}[Proof of \cref{W2contraction}]\label{proofw2contraction}
    Concerning the moments, we have
    \begin{align}
        \bigg\langle\frac{d}{dt}\mu^{X}_{t},\|\cdot\|^2_{2}\bigg\rangle&=\bigg\langle\mu^{X}_{t},\mathcal{L}_{\mu^{X}_{t}}\|\cdot\|^2_{2}\bigg\rangle\\
        &=\int\bigg(2d-2\nabla\frac{\delta H}{\delta m}(\mu^{X}_{t},x)\cdot x\bigg)\mu^{X}_{t}(dx)\nonumber\\
        &\leq 2d-2\bigg[\bigg(\lambda-\bigg\|\nabla\frac{\delta H}{\delta m}(\mu^{X}_{t},0)\bigg\|_{2}\bigg)\bigg\langle\mu^{X}_{t},\|\cdot\|^2_{2}\bigg\rangle-\bigg\|\nabla\frac{\delta H}{\delta m}(\mu^{X}_{t},0)\bigg\|_{2}\bigg]\nonumber\\
        &\leq 2d+2\bigg\|\nabla\frac{\delta H}{\delta m}(\mu^{X}_{t},0)\bigg\|_{2}-2\bigg(\lambda-\bigg\|\nabla\frac{\delta H}{\delta m}(\mu^{X}_{t},0)\bigg\|_{2}\bigg)\bigg\langle\mu^{X}_{t},\|\cdot\|^2_{2}\bigg\rangle.\nonumber
    \end{align}
    We deduce that there exists $(\alpha,\beta)\in\mathbb{R}^{2}$ such that
    \begin{equation}
        \bigg\langle\frac{d}{dt}\mu^{X}_{t},\|\cdot\|^2_{2}\bigg\rangle\leq\alpha+\beta\bigg\langle\mu^{X}_{t},\|\cdot\|^2_{2}\bigg\rangle.
    \end{equation}
    By Gronwall's lemma, we have
    \begin{equation}\label{bounded2momentum}
        \forall t\geq0,\quad\bigg\langle\mu^{X}_{t},\|\cdot\|^2_{2}\bigg\rangle\leq\begin{cases}
            \bigg\langle\mu^{X}_{0},\|\cdot\|^2_{2}\bigg\rangle e^{\beta t}+\frac{\alpha}{\beta}(e^{\beta t}-1)\quad&\textnormal{if}\quad\beta\neq0;\\
            \bigg\langle\mu^{X}_{0},\|\cdot\|^2_{2}\bigg\rangle+\alpha t\quad&\textnormal{if}\quad\beta=0.
        \end{cases}
    \end{equation}
    As for the exponential contraction, to do this, we consider the coupling given by
    \begin{equation}
        \begin{cases}
            dX_{t}&=-\nabla\frac{\delta H}{\delta m}(\mu^{X}_{t},X_{t})dt+\sqrt{2}dB_{t}\\
            dY_{t}&=-\nabla\frac{\delta H}{\delta m}(\mu^{Y}_{t},Y_{t})dt+\sqrt{2}dB_{t}.
        \end{cases}
    \end{equation}
    We have
    \begin{align}
        \frac{1}{2}\frac{d}{dt}\|X_{t}-Y_{t}\|^2_{2}&=\bigg\langle\frac{d}{dt}(X_{t}-Y_{t}),X_{t}-Y_{t}\bigg\rangle\\
        &=-\bigg\langle\nabla\frac{\delta H}{\delta m}(\mu^{X}_{t},X_{t})-\nabla\frac{\delta H}{\delta m}(\mu^{Y}_{t},Y_{t}),X_{t}-Y_{t}\bigg\rangle\nonumber\\
        &=-\bigg\langle\nabla\frac{\delta H}{\delta m}(\mu^{Y}_{t},X_{t})-\nabla\frac{\delta H}{\delta m}(\mu^{Y}_{t},Y_{t}),X_{t}-Y_{t}\bigg\rangle\nonumber\\
        &-\bigg\langle\nabla\frac{\delta H}{\delta m}(\mu^{X}_{t},X_{t})-\nabla\frac{\delta H}{\delta m}(\mu^{Y}_{t},X_{t}),X_{t}-Y_{t}\bigg\rangle.\nonumber
    \end{align}
    On the one hand, according to $\mathbf{(i)}$, we have
    \begin{equation}
        -\bigg\langle\nabla\frac{\delta H}{\delta m}(\mu^{Y}_{t},X_{t})-\nabla\frac{\delta H}{\delta m}(\mu^{Y}_{t},Y_{t}),X_{t}-Y_{t}\bigg\rangle\leq-\lambda\|X_{t}-Y_{t}\|^2_{2}.
    \end{equation}
    On the other hand, by definition of the flat derivative and according to $\mathbf{(ii)}$, we have
    \begin{align}
       & -\bigg\langle\nabla\frac{\delta H}{\delta m}(\mu^{X}_{t},X_{t})-\nabla\frac{\delta H}{\delta m}(\mu^{Y}_{t},X_{t}),X_{t}-Y_{t}\bigg\rangle\\
       &=-\int_{0}^{1}\int\bigg\langle\nabla_{x}\frac{\delta^2H}{\delta m^2}(s\mu^{X}_{t}+(1-s)\mu^{Y}_{t},X_{t},y),X_{t}-Y_{t}\bigg\rangle(\mu^{X}_{t}-\mu^{Y}_{t})(dy)ds\nonumber\\
        &\leq\|\mathcal{D}^{2}_{m}H\|_{\mathbf{op},\infty}\|X_{t}-Y_{t}\|_{2}\mathcal{W}_{1}(\mu^{X}_{t},\mu^{Y}_{t})\nonumber\\
        &\leq\|\mathcal{D}^{2}_{m}H\|_{\mathbf{op},\infty}\|X_{t}-Y_{t}\|_{2}\mathbb{E}[\|X_{t}-Y_{t}\|_{2}].\nonumber
    \end{align}
We deduce that
\begin{align}
    \frac{1}{2}\frac{d}{dt}\|X_{t}-Y_{t}\|^2_{2}&\leq-\lambda\|X_{t}-Y_{t}\|^2_{2}+\|\mathcal{D}^{2}_{m}H\|_{\mathbf{op},\infty}\|X_{t}-Y_{t}\|_{2}\mathbb{E}[\|X_{t}-Y_{t}\|_{2}];\\
    \frac{1}{2}\frac{d}{dt}\mathbb{E}[\|X_{t}-Y_{t}\|^2_{2}]&\leq-\lambda\mathbb{E}[\|X_{t}-Y_{t}\|^2_{2}]+\|\mathcal{D}^{2}_{m}H\|_{\mathbf{op},\infty}\bigg(\mathbb{E}[\|X_{t}-Y_{t}\|_{2}]\bigg)^2\\
    &\leq (\|\mathcal{D}^{2}_{m}H\|_{\mathbf{op},\infty}-\lambda)\mathbb{E}[\|X_{t}-Y_{t}\|^2_{2}]\quad\textnormal{Jensen's inequality}.\nonumber
\end{align}
According to Gronwall's lemma, we have
\begin{equation}
    \forall t\geq0,\quad\mathbb{E}[\|X_{t}-Y_{t}\|^2_{2}]\leq\mathbb{E}[\|X_{0}-Y_{0}\|^2_{2}]e^{2(\|\mathcal{D}^{2}_{m}H\|_{\mathbf{op},\infty}-\lambda)t}.
\end{equation}
We deduce that
\begin{equation}
    \forall t\geq0,\quad\mathcal{W}^2_{2}(\mu^{X}_{t},\mu^{Y}_{t})\leq\mathbb{E}[\|X_{t}-Y_{t}\|^2_{2}]\leq\mathbb{E}[\|X_{0}-Y_{0}\|^2_{2}]e^{2(\|\mathcal{D}^{2}_{m}H\|_{\mathbf{op},\infty}-\lambda)t}.
\end{equation}
We conclude that
\begin{equation}
    \forall t\geq0,\quad\mathcal{W}_{2}(\mu^{X}_{t},\mu^{Y}_{t})\leq e^{(\|\mathcal{D}^{2}_{m}H\|_{\mathbf{op},\infty}-\lambda)t}\mathcal{W}_{2}(\mu^{X}_{0},\mu^{Y}_{0}).
\end{equation}
\end{proof}

\begin{proof}[Proof of \cref{uniformpropachaosmeanfieldLangevin}]\label{proofuniformpropachaosmeanfieldLangevin}
    Consider $(\widetilde{X}^{(n)})_{n\geq1}$ a sequence of independent copies of the mean field Langevin process $(X_{t})_{t\geq0}$. By synchronous coupling of the particle system with these copies, we have
    \begin{align}
        \frac{d}{dt}\sum_{i=1}^{n}\|X^{(i)}_{t}-\widetilde{X}^{(i)}_{t}\|^2&=-2\sum_{i=1}^{n}\bigg\langle\nabla\frac{\delta H}{\delta m}(\mu_{\mathbf{X}_{t}},X^{(i)}_{t})-\nabla\frac{\delta H}{\delta m}(\mu^{X}_{t},\widetilde{X}^{(i)}_{t}),X^{(i)}_{t}-\widetilde{X}^{(i)}_{t}\bigg\rangle\\
        &\leq-2\lambda\sum_{i=1}^{n}\|X^{(i)}_{t}-\widetilde{X}^{(i)}_{t}\|^2-2\sum_{i=1}^{n}\bigg\langle\nabla\frac{\delta H}{\delta m}(\mu_{\mathbf{X}_{t}},\widetilde{X}^{(i)}_{t})-\nabla\frac{\delta H}{\delta m}(\mu^{X}_{t},\widetilde{X}^{(i)}_{t}),X^{(i)}_{t}-\widetilde{X}^{(i)}_{t}\bigg\rangle\nonumber\\
        &\leq-2\lambda\sum_{i=1}^{n}\|X^{(i)}_{t}-\widetilde{X}^{(i)}_{t}\|^2+2\|\mathcal{D}^2_{m}H\|_{\mathbf{op},\infty}\mathcal{W}_{1}(\mu_{\mathbf{X}_{t}},\mu^{X}_{t})\sum_{i=1}^{n}\|X^{(i)}_{t}-\widetilde{X}^{(i)}_{t}\|.\nonumber
    \end{align}
    On the other hand, we have
    \begin{align}
        \mathcal{W}_{1}(\mu_{\mathbf{X}_{t}},\mu^{X}_{t})&\leq\mathcal{W}_{2}(\mu_{\mathbf{X}_{t}},\mu^{X}_{t});\\
        \mathcal{W}_{2}(\mu_{\mathbf{X}_{t}},\mu^{X}_{t})&\leq\mathcal{W}_{2}(\mu_{\mathbf{X}_{t}},\mu_{\widetilde{\mathbf{X}}_{t}})+\mathcal{W}_{2}(\mu_{\widetilde{\mathbf{X}}_{t}},\mu^{X}_{t}).
    \end{align}
    By Young's inequality, we have
    \begin{align}
        \frac{d}{dt}\sum_{i=1}^{n}\|X^{(i)}_{t}-\widetilde{X}^{(i)}_{t}\|^2&\leq2(\|\mathcal{D}^2_{m}H\|_{\mathbf{op},\infty}-\lambda)\sum_{i=1}^{n}\|X^{(i)}_{t}-\widetilde{X}^{(i)}_{t}\|^2\\
        &+n\|\mathcal{D}^2_{m}H\|_{\mathbf{op},\infty}\bigg(\mathcal{W}^2_{2}(\mu_{\mathbf{X}_{t}},\mu_{\widetilde{\mathbf{X}}_{t}})+\mathcal{W}^2_{2}(\mu_{\widetilde{\mathbf{X}}_{t}},\mu^{X}_{t})\bigg).\nonumber
    \end{align}
    We have
    \begin{equation}
        \mathbb{E}[\mathcal{W}^2_{2}(\mu_{\mathbf{X}_{t}},\mu_{\widetilde{\mathbf{X}}_{t}})]\leq\frac{1}{n}\sum_{i=1}^{n}\mathbb{E}[\|X^{(i)}_{t}-\widetilde{X}^{(i)}_{t}\|^2].
    \end{equation}
    By symmetry, we have
    \begin{equation}
        \forall j\in\{1,\ldots,n\},\quad\mathbb{E}[\|X^{(j)}_{t}-\widetilde{X}^{(j)}_{t}\|^2]=\frac{1}{n}\sum_{i=1}^{n}\mathbb{E}[\|X^{(i)}_{t}-\widetilde{X}^{(i)}_{t}\|^2].
    \end{equation}
    Setting $\gamma_{t}:=\mathbb{E}[\|X^{(1)}_{t}-\widetilde{X}^{(1)}_{t}\|^2]$, we have
    \begin{equation}
        \frac{d}{dt}\gamma_{t}\leq(3\|\mathcal{D}^2_{m}H\|_{\mathbf{op},\infty}-2\lambda)\gamma_{t}+\|\mathcal{D}^2_{m}H\|_{\mathbf{op},\infty}\mathbb{E}[\mathcal{W}^2_{2}(\mu_{\widetilde{\mathbf{X}}_{t}},\mu^{X}_{t})].
    \end{equation}
    By \cref{bounded2momentum}, we have
    \begin{align}
        \forall t\geq0,\quad\bigg\langle\mu^{X}_{t},\|\cdot\|^2_{2}\bigg\rangle&\leq\begin{cases}
            \bigg\langle\mu^{X}_{0},\|\cdot\|^2_{2}\bigg\rangle e^{\beta t}+\frac{\alpha}{\beta}(e^{\beta t}-1)\quad&\textnormal{if}\quad\beta\neq0;\\
            \bigg\langle\mu^{X}_{0},\|\cdot\|^2_{2}\bigg\rangle+\alpha t\quad&\textnormal{if}\quad\beta=0.
        \end{cases}\\
        \alpha&:=2d+2\sup_{\mu}\bigg\|\nabla\frac{\delta H}{\delta m}(\mu,0)\bigg\|_{2};\\
        \beta&:=-2\bigg(\lambda-\sup_{\mu}\bigg\|\nabla\frac{\delta H}{\delta m}(\mu,0)\bigg\|_{2}\bigg).
    \end{align}
    If $\beta<0$, we have
    \begin{equation}
        \sup_{t\geq0}\bigg\langle\mu^{X}_{t},\|\cdot\|^2_{2}\bigg\rangle\leq\bigg\langle\mu^{X}_{0},\|\cdot\|^2_{2}\bigg\rangle-\frac{\alpha}{\beta}.
    \end{equation}
     By the result of Fournier and Guillin (\cite[Theorem 1]{fournier2013rateconvergencewassersteindistance}), we have
     \begin{align}
         \mathbb{E}[\mathcal{W}^2_{2}(\mu_{\widetilde{\mathbf{X}}_{t}},\mu^{X}_{t})]&\leq C(d)\bigg(\bigg\langle\mu^{X}_{0},\|\cdot\|^2_{2}\bigg\rangle-\frac{\alpha}{\beta}\bigg)\delta_{d}(n);\\
         \delta_{d}(n)&:=\begin{cases}
             \frac{1}{\sqrt{n}}\quad&\textnormal{if}\quad d<4;\\
             \frac{\ln(n+1)}{\sqrt{n}}\quad&\textnormal{if}\quad d=4;\\
             \frac{1}{n^{\frac{2}{d}}}\quad&\textnormal{if}\quad d>4.
         \end{cases}
     \end{align}
     By Gronwall's lemma, we have ($\gamma_{0}=0$)
     \begin{align}
        \forall t\geq0,\quad\gamma_{t}&\leq\begin{cases}
            \frac{\alpha_{1}}{\beta_{1}}(e^{\beta_{1} t}-1)\quad&\textnormal{if}\quad\beta_{1}\neq0;\\
            \alpha_{1} t\quad&\textnormal{if}\quad\beta_{1}=0.
        \end{cases}\\
        \alpha_{1}&:=\|\mathcal{D}^2_{m}H\|_{\mathbf{op},\infty}C(d)\bigg(\bigg\langle\mu^{X}_{0},\|\cdot\|^2_{2}\bigg\rangle-\frac{\alpha}{\beta}\bigg)\delta_{d}(n);\\
        \beta_{1}&:=3\|\mathcal{D}^2_{m}H\|_{\mathbf{op},\infty}-2\lambda.
    \end{align}
    In this case, the bound is uniform in time if, and only if, $\beta_{1}<0$ :
    \begin{equation}
        \sup_{t\geq0}\gamma_{t}\leq-\frac{\alpha_{1}}{\beta_{1}}.
    \end{equation}

    If $\beta=0$, we have:
    \begin{equation}
        \sup_{t\in [0,T]}\bigg\langle\mu^{X}_{t},\|\cdot\|^2_{2}\bigg\rangle\leq\bigg\langle\mu^{X}_{0},\|\cdot\|^2_{2}\bigg\rangle+\alpha T.
    \end{equation}
We deduce that
\begin{align}
        \forall t\in [0,T],\quad\gamma_{t}&\leq\begin{cases}
            \frac{\alpha_{2}}{\beta_{1}}(e^{\beta_{1} t}-1)\quad&\textnormal{if}\quad\beta_{1}\neq0;\\
            \alpha_{2} t\quad&\textnormal{if}\quad\beta_{1}=0.
        \end{cases}\\
        \alpha_{2}&:=\|\mathcal{D}^2_{m}H\|_{\mathbf{op},\infty}C(d)\bigg(\bigg\langle\mu^{X}_{0},\|\cdot\|^2_{2}\bigg\rangle+\alpha T\bigg)\delta_{d}(n);\\
        \beta_{1}&:=3\|\mathcal{D}^2_{m}H\|_{\mathbf{op},\infty}-2\lambda.
    \end{align}
    In this case, since the moments are only bounded in finite time horizon, the propagation of chaos cannot be uniform in time.\newline

    If $\beta>0$, we have
    \begin{equation}
        \sup_{t\in [0,T]}\bigg\langle\mu^{X}_{t},\|\cdot\|^2_{2}\bigg\rangle\leq\langle\mu^{X}_{0},\|\cdot\|^2_{2}\bigg\rangle e^{\beta T}+\frac{\alpha}{\beta}(e^{\beta T}-1).
    \end{equation}
We deduce that
\begin{align}
        \forall t\in [0,T],\quad\gamma_{t}&\leq\begin{cases}
            \frac{\alpha_{3}}{\beta_{1}}(e^{\beta_{1} t}-1)\quad&\textnormal{if}\quad\beta_{1}\neq0;\\
            \alpha_{3} t\quad&\textnormal{if}\quad\beta_{1}=0.
        \end{cases}\\
        \alpha_{3}&:=\|\mathcal{D}^2_{m}H\|_{\mathbf{op},\infty}C(d)\bigg(\bigg\langle\mu^{X}_{0},\|\cdot\|^2_{2}\bigg\rangle e^{\beta T}+\frac{\alpha}{\beta}(e^{\beta T}-1)\bigg)\delta_{d}(n);\\
        \beta_{1}&:=3\|\mathcal{D}^2_{m}H\|_{\mathbf{op},\infty}-2\lambda.
    \end{align}
    In this case, since the moments are only bounded in finite time horizon, the propagation of chaos cannot be uniform in time.
    
\end{proof}

\subsection{Proof of main results of the kinetic case}
\begin{proof}[Proof of \cref{kineticw2contraction}]\label{proofkineticw2contraction}
    Let $(P_{t},V_{t})_{t\geq0}$ and $(\overline{P}_{t},\overline{V}_{t})_{t\geq0}$ two solutions of
    \begin{equation}
        \begin{cases}
            dP_{t}&=V_{t}dt;\\
            dV_{t}&=\sqrt{2}dB_{t}-A(V_{t})dt-B(P_{t})dt-\nabla\frac{\delta H}{\delta m}(\mathbb{P}_{P_{t}},P_{t})dt;
        \end{cases}
    \end{equation}
    such that $\mathbb{P}_{(P_{0},V_{0})}\in\mathcal{P}_{2}(\mathbb{R}^{d}\times\mathbb{R}^{d})$ and $\mathbb{P}_{(\overline{P}_{0},\overline{V}_{0})}\in\mathcal{P}_{2}(\mathbb{R}^{d}\times\mathbb{R}^{d})$. Then, by difference, $(p_{t},v_{t}):=(P_{t}-\overline{P}_{t},V_{t}-\overline{V}_{t})$ evolves according to
    \begin{equation}
        \begin{cases}
            dp_{t}&=v_{t}dt;\\
            dv_{t}&=-\bigg(A(V_{t})-A(\overline{V}_{t})+\lambda_{B}p_{t}+D(P_{t})-D(\overline{P}_{t})+\nabla\frac{\delta H}{\delta m}(\mathbb{P}_{P_{t}},P_{t})-\nabla\frac{\delta H}{\delta m}(\mathbb{P}_{\overline{P}_{t}},\overline{P}_{t})\bigg)dt.
        \end{cases}
    \end{equation}
    Then, if $a$ and $b$ are positive constants to be chosen later on,
    \begin{align}
        \frac{d}{dt}\bigg(a|p_{t}|^2+2p_{t}\cdot v_{t}+b|v_{t}|^2\bigg)&=2ap_{t}\cdot v_{t}+2|v_{t}|^2\\
        &-2p_{t}\cdot\bigg(A(V_{t})-A(\overline{V}_{t})+\lambda_{B}p_{t}+D(P_{t})-D(\overline{P}_{t})\bigg)\nonumber\\
        &-2bv_{t}\cdot\bigg(A(V_{t})-A(\overline{V}_{t})+\lambda_{B}p_{t}+D(P_{t})-D(\overline{P}_{t})\bigg)\nonumber\\
        &-2(p_{t}+bv_{t})\cdot\bigg(\nabla\frac{\delta H}{\delta m}(\mathbb{P}_{P_{t}},P_{t})-\nabla\frac{\delta H}{\delta m}(\mathbb{P}_{\overline{P}_{t}},\overline{P}_{t})\bigg).\nonumber
    \end{align}
    By the Cauchy-Schwarz inequality and assumptions on $A$ and $D$, the first four terms are bounded by above by
    \begin{equation}
        2(a-\lambda_{B}b)p_{t}\cdot v_{t}+2([A]_{1}+[D]_{1}b)|p_{t}|\cdot |v_{t}|-2(\lambda_{A}b-1)|v_{t}|^2-2(\lambda_{B}-[D]_{1})|p_{t}|^2.
    \end{equation}
    As for the last remaining term, we have
    \begin{align}
        -2(p_{t}+bv_{t})\cdot\bigg(\nabla\frac{\delta H}{\delta m}(\mathbb{P}_{P_{t}},P_{t})-\nabla\frac{\delta H}{\delta m}(\mathbb{P}_{\overline{P}_{t}},\overline{P}_{t})\bigg)&=-2(p_{t}+bv_{t})\cdot\bigg(\nabla\frac{\delta H}{\delta m}(\mathbb{P}_{P_{t}},P_{t})-\nabla\frac{\delta H}{\delta m}(\mathbb{P}_{{P}_{t}},\overline{P}_{t})\bigg)\\
        &-2(p_{t}+bv_{t})\cdot\bigg(\nabla\frac{\delta H}{\delta m}(\mathbb{P}_{P_{t}},\overline{P}_{t})-\nabla\frac{\delta H}{\delta m}(\mathbb{P}_{\overline{P}_{t}},\overline{P}_{t})\bigg).\nonumber
    \end{align}
    Then, we have
    \begin{align}
        -2p_{t}\cdot\bigg(\nabla\frac{\delta H}{\delta m}(\mathbb{P}_{P_{t}},P_{t})-\nabla\frac{\delta H}{\delta m}(\mathbb{P}_{{P}_{t}},\overline{P}_{t})\bigg)&\leq2[\mathcal{D}_{m}H]_{1,\infty}|p_{t}|^2;\\
        -2p_{t}\cdot\bigg(\nabla\frac{\delta H}{\delta m}(\mathbb{P}_{P_{t}},\overline{P}_{t})-\nabla\frac{\delta H}{\delta m}(\mathbb{P}_{\overline{P}_{t}},\overline{P}_{t})\bigg)&\leq2\|\mathcal{D}^{2}_{m}H\|_{\mathbf{op},\infty}|p_{t}|\mathbb{E}[|p_{t}|];\\
        -2bv_{t}\cdot\bigg(\nabla\frac{\delta H}{\delta m}(\mathbb{P}_{P_{t}},P_{t})-\nabla\frac{\delta H}{\delta m}(\mathbb{P}_{{P}_{t}},\overline{P}_{t})\bigg)&\leq2b[\mathcal{D}_{m}H]_{1,\infty}|p_{t}|\cdot |v_{t}|;\\
        -2bv_{t}\cdot\bigg(\nabla\frac{\delta H}{\delta m}(\mathbb{P}_{P_{t}},\overline{P}_{t})-\nabla\frac{\delta H}{\delta m}(\mathbb{P}_{\overline{P}_{t}},\overline{P}_{t})\bigg)&\leq2b\|\mathcal{D}^{2}_{m}H\|_{\mathbf{op},\infty}|v_{t}|\mathbb{E}[|p_{t}|].
    \end{align}
    By the identity $2ab\leq a^2+b^2$, we have
    \begin{equation}
        -2bv_{t}\cdot\bigg(\nabla\frac{\delta H}{\delta m}(\mathbb{P}_{P_{t}},P_{t})-\nabla\frac{\delta H}{\delta m}(\mathbb{P}_{{P}_{t}},\overline{P}_{t})\bigg)\leq b[\mathcal{D}_{m}H]_{1,\infty}(|p_{t}|^2+|v_{t}|^2).
    \end{equation}
    By Jensen's inequality and Young's inequality, we have
    \begin{align}
        -2\mathbb{E}\bigg[p_{t}\cdot\bigg(\nabla\frac{\delta H}{\delta m}(\mathbb{P}_{P_{t}},\overline{P}_{t})-\nabla\frac{\delta H}{\delta m}(\mathbb{P}_{\overline{P}_{t}},\overline{P}_{t})\bigg)\bigg]&\leq2\|\mathcal{D}^{2}_{m}H\|_{\mathbf{op},\infty}\mathbb{E}[|p_{t}|^{2}];\\
        -2b\mathbb{E}\bigg[v_{t}\cdot\bigg(\nabla\frac{\delta H}{\delta m}(\mathbb{P}_{P_{t}},\overline{P}_{t})-\nabla\frac{\delta H}{\delta m}(\mathbb{P}_{\overline{P}_{t}},\overline{P}_{t})\bigg)\bigg]&\leq b\|\mathcal{D}^{2}_{m}H\|_{\mathbf{op},\infty}\mathbb{E}[|p_{t}|^{2}+|v_{t}|^{2}].
    \end{align}
  Collecting all terms leads to the bound
  \begin{align}
      \frac{d}{dt}\mathbb{E}\bigg[a|p_{t}|^2+2p_{t}\cdot v_{t}+b|v_{t}|^2\bigg]&\leq2\bigg(a-\lambda_{B}b\bigg)\mathbb{E}[p_{t}\cdot v_{t}]+2\bigg([A]_{1}+[D]_{1}b\bigg)\mathbb{E}[|p_{t}|\cdot|v_{t}|]\\
      &-2\bigg(\lambda_{B}-[D]_{1}-([\mathcal{D}_{m}H]_{1,\infty}+\|\mathcal{D}^{2}_{m}H\|_{\mathbf{op},\infty})\nonumber\\
      &-\frac{b}{2}([\mathcal{D}_{m}H]_{1,\infty}+\|\mathcal{D}^{2}_{m}H\|_{\mathbf{op},\infty})\bigg)\mathbb{E}[|p_{t}|^2]\nonumber\\
      &-2\bigg(\lambda_{A}b-1-\frac{b}{2}([\mathcal{D}_{m}H]_{1,\infty}+\|\mathcal{D}^{2}_{m}H\|_{\mathbf{op},\infty})\bigg)\mathbb{E}[|v_{t}|^2]\nonumber.
  \end{align}
   For all $a\geq\lambda_{B}b$, by the Cauchy-Schwarz inequality and the Young's inequality ($2ab\leq a^2+b^2$), we have
    \begin{align}
        \frac{d}{dt}\mathbb{E}\bigg[a|p_{t}|^2+2p_{t}\cdot v_{t}+b|v_{t}|^2\bigg]&\leq\bigg(a-\lambda_{B}b+[A]_{1}+[D]_{1}b\\
        &-2\bigg(\lambda_{B}-[D]_{1}-([\mathcal{D}_{m}H]_{1,\infty}+\|\mathcal{D}^{2}_{m}H\|_{\mathbf{op},\infty})\nonumber\\
        &-\frac{b}{2}([\mathcal{D}_{m}H]_{1,\infty}+\|\mathcal{D}^{2}_{m}H\|_{\mathbf{op},\infty})\bigg)\bigg)\mathbb{E}[|p_{t}|^2]\nonumber\\
        &+\bigg(a-\lambda_{B}b+[A]_{1}+[D]_{1}b\nonumber\\
        &-2\bigg(\lambda_{A}b-1-\frac{b}{2}([\mathcal{D}_{m}H]_{1,\infty}+\|\mathcal{D}^{2}_{m}H\|_{\mathbf{op},\infty})\bigg)\bigg)\mathbb{E}[|v_{t}|^2]\nonumber.
    \end{align}
    For all $b>\frac{1}{\sqrt{\lambda_{B}}}$, $Q_{\lambda_{B}b,b}: (p,v)\in\mathbb{R}^{2d}\mapsto\lambda_{B}b|p|^2+2p\cdot v+b|v|^2$ is a positive quadratic
    form on $\mathbb{R}^{2d}$. Then we look for $a$ and $b$ such that
    \begin{equation}
        \begin{cases}
            a-\lambda_{B}b+[A]_{1}+[D]_{1}b-2\bigg(\lambda_{B}-[D]_{1}-([\mathcal{D}_{m}H]_{1,\infty}+\|\mathcal{D}^{2}_{m}H\|_{\mathbf{op},\infty})-\frac{b}{2}([\mathcal{D}_{m}H]_{1,\infty}+\|\mathcal{D}^{2}_{m}H\|_{\mathbf{op},\infty})\bigg)&<0;\\
            a-\lambda_{B}b+[A]_{1}+[D]_{1}b-2\bigg(\lambda_{A}b-1-\frac{b}{2}([\mathcal{D}_{m}H]_{1,\infty}+\|\mathcal{D}^{2}_{m}H\|_{\mathbf{op},\infty})\bigg)&<0;\\
            a\geq\lambda_{B}b.
        \end{cases}
    \end{equation}
    in such a way that
    \begin{equation}
        \exists C>0\quad\forall t\geq0,\quad\frac{d}{dt}\mathbb{E}\bigg[Q_{a,b}\bigg((p_{t},v_{t})\bigg)\bigg]\leq -C\mathbb{E}[|p_{t}|^2+|v_{t}|^2].
    \end{equation}
    Let us set $\gamma:=[\mathcal{D}_{m}H]_{1,\infty}+\|\mathcal{D}^{2}_{m}H\|_{\mathbf{op},\infty}$.
    If $a=\lambda_{B}b$,  $b$ satisfies
    \begin{equation}
        \begin{cases}
            b<\frac{2\lambda_{B}-2[D]_{1}-[A]_{1}-2\gamma}{[D]_{1}+\gamma};\\
            ([D]_{1}+\gamma-2\lambda_{A})b<-[A]_{1}-2.
        \end{cases}
    \end{equation}
    Such $b$ exists given that
    \begin{align}
        2\lambda_{B}-2[D]_{1}-[A]_{1}-2\gamma&>0;\\
        [D]_{1}+\gamma-2\lambda_{A}&<0;\\
        \frac{-[A]_{1}-2}{[D]_{1}+\gamma-2\lambda_{A}}&<\frac{2\lambda_{B}-2[D]_{1}-[A]_{1}-2\gamma}{[D]_{1}+\gamma}.
    \end{align}
    In fact, for every positive $\varepsilon$, we have
    \begin{equation}
        \frac{d}{dt}\mathbb{E}\bigg[\lambda_{B}b|p_{t}|^2+2p_{t}\cdot v_{t}+b|v_{t}|^2\bigg]\leq-\bigg(2\lambda_{B}-[A]_{1}-2\eta-\varepsilon-\eta b\bigg)\mathbb{E}[|p_{t}|^{2}]-\bigg((2\lambda_{A}-\eta)b-2-[A]_{1}-\frac{[A]^{2}_{1}}{\varepsilon}\bigg)\mathbb{E}[|v_{t}|^{2}]
    \end{equation}
    where $\eta:=\gamma+[D]_{1}$. For us to have
    \begin{equation}
        \begin{cases}
            2\lambda_{B}-[A]_{1}-2\eta-\varepsilon-\eta b&>0;\\
            (2\lambda_{A}-\eta)b-2-[A]_{1}-\frac{[A]^{2}_{1}}{\varepsilon}&>0;
        \end{cases}
    \end{equation}
    it is necessary that $\eta<2\lambda_{A}$, which is assumed in the sequel.  Then, for instance for $\varepsilon=\lambda_{B}$, the conditions are equivalent to
    \begin{equation}
        \frac{2+[A]_{1}+\frac{[A]^{2}_{1}}{\lambda_{B}}}{2\lambda_{A}-\eta}<b<\frac{\lambda_{B}-[A]_{1}-2\eta}{\eta}\quad\textnormal{and}\quad\eta<2\lambda_{A}.
    \end{equation}
    We look for $\eta$ such that
    \begin{equation}
        \frac{2+[A]_{1}+\frac{[A]^{2}_{1}}{\lambda_{B}}}{2\lambda_{A}-\eta}<\frac{\lambda_{B}-[A]_{1}-2\eta}{\eta},
    \end{equation}
    that is,
    \begin{equation}
        2\eta^2-\eta\bigg(2+\frac{[A]^{2}_{1}}{\lambda_{B}}+\lambda_{B}+4\lambda_{A}\bigg)+2\lambda_{A}(\lambda_{B}-[A]_{1})>0.
    \end{equation}
    This polynomial takes negative values at $\eta=2\lambda_{A}$, so it is positive on an interval $[0,\eta_{0}[$ for some $\eta_{0}<2\lambda_{A}$: there exists $b$ with all the above conditions for any $\eta\in [0,\eta_{0}[$.
    Hence there exists a constant $\eta_{0}$, depending only on $[A]_{1}$, $\lambda_{A}$ and $\lambda_{B}$, such that, if $\gamma+[D]_{1}<\eta_{0}$, then there exist a positive quadratic form $Q_{\lambda_{B}b,b}$ on $\mathbb{R}^{2d}$ and a constant $C$ depending only on $[A]_{1}$, $\lambda_{A}$, $\lambda_{B}$, $\gamma$ and $[D]_{1}$ such that
    \begin{equation}
        \forall t\geq0,\quad\frac{d}{dt}\mathbb{E}\bigg[Q_{\lambda_{B}b,b}\bigg((p_{t},v_{t})\bigg)\bigg]\leq -C\mathbb{E}[|p_{t}|^2+|v_{t}|^2].
    \end{equation}
    In turn, since $\sqrt{Q_{\lambda_{B}b,b}}$ and $(p,v)\mapsto \sqrt{|p|^2+|v|^2}$ are equivalent on $\mathbb{R}^{2d}$,  this is bounded by $-2C\mathbb{E}[Q_{\lambda_{B}b,b}((p_{t},v_{t}))]$ for a new constant $C>0$, so that,
    \begin{equation}
        \forall t\geq0,\quad \mathbb{E}\bigg[Q_{\lambda_{B}b,b}\bigg((P_{t},V_{t})-(\overline{P}_{t},\overline{V}_{t})\bigg)\bigg]\leq e^{-2Ct}\mathbb{E}\bigg[Q_{\lambda_{B}b,b}\bigg((P_{0},V_{0})-(\overline{P}_{0},\overline{V}_{0})\bigg)\bigg]
    \end{equation}
    by integration. We deduce by equivalence of norms on the space of configurations and the existence of an optimal coupling that
    \begin{equation}
        \forall t\geq0,\quad\mathcal{W}_{2}(\mathbb{P}_{(P_{t},V_{t})},\mathbb{P}_{(\overline{P}_{t},\overline{V}_{t})})\leq C'e^{-Ct}\mathcal{W}_{2}(\mathbb{P}_{(P_{0},V_{0})},\mathbb{P}_{(\overline{P}_{0},\overline{V}_{0})}).
    \end{equation}
    We conclude that there exists a unique equilibrium (\cite[Lemma 10]{BolleyGuillinMalrieu}) and that we contract exponentially towards the latter.
\end{proof}

\begin{proof}[Proof of \cref{uniformpropachaoskineticmeanfieldLangevin}]\label{proofuniformpropachaoskineticmeanfieldLangevin}
The time-uniform propagation of chaos requires a time-uniform bound of the second
moment of $\mathbb{P}_{(P_{t},V_{t})}$. Let $(f_{t})_{t\geq0}$ be a solution to \cref{meanfieldkineticPDE} with initial datum $f_{0}\in\mathcal{P}_{2}(\mathbb{R}^{d}\times\mathbb{R}^{d})$ and let $a$ and $b$ be positive numbers to be chosen later on. Note that for all test function $\varphi$,
\begin{align}
    &\frac{d}{dt}\int\varphi(p,v)f_{t}(dpdv)\\
    &=\int\bigg(\Delta_{v}\varphi(p,v)+v\cdot\nabla_{p}\varphi(p,v)-\bigg(A(v)+B(p)+\mathcal{D}_{m}H(f_{t},p)\bigg)\cdot\nabla_{v}\varphi(p,v)\bigg)f_{t}(dpdv).\nonumber
\end{align}
Then
\begin{align}
    &\frac{d}{dt}\int\bigg(a|p|^2+2p\cdot v+b|v|^2\bigg)f_{t}(dpdv)\\
    &=2bd+2\int\bigg( v\cdot(ap+v)-(p+bv)\cdot\bigg(A(v)+B(p)+\mathcal{D}_{m}H(f_{t},p)\bigg)\bigg)f_{t}(dpdv)\nonumber
\end{align}
where, by Young’s inequality and assumption on $A$, $B$ and $H$, 
\begin{align}
    -2p\cdot A(v)&=-2p\cdot (A(v)-A(0))-2p\cdot A(0)\\
    &\leq2[A]_{1}|p|\cdot|v|-2p\cdot A(0);\nonumber\\
    -2p\cdot B(p)&=-2p\cdot(\lambda_{B}p+D(p)-D(0)+D(0))\\
    &\leq-(2\lambda_{B}-2[D]_{1})|p|^2-2p\cdot D(0);\nonumber\\
    -2bv\cdot A(v)&=-2bv\cdot (A(v)-A(0)+A(0))\\
    &\leq-2b\lambda_{A}|v|^2-2bv\cdot A(0);\nonumber\\
    -2bv\cdot B(p)&=-2bv\cdot(\lambda_{B}p+D(p)-D(0)+D(0))\\
    &\leq-2b\lambda_{B}p\cdot v+2b[D]_{1}|p|\cdot|v|-2bv\cdot D(0);\nonumber\\
    -2p\cdot\mathcal{D}_{m}H(f_{t},p)&=-2p\cdot(\mathcal{D}_{m}H(f_{t},p)-\mathcal{D}_{m}H(f_{t},0)+\mathcal{D}_{m}H(f_{t},0))\\
    &\leq 2[\mathcal{D}_{m}H]_{1,\infty}|p|^2-2p\cdot\mathcal{D}_{m}H(f_{t},0);\nonumber\\
    -2bv\cdot\mathcal{D}_{m}H(f_{t},p)&=-2bv\cdot(\mathcal{D}_{m}H(f_{t},p)-\mathcal{D}_{m}H(f_{t},0)+\mathcal{D}_{m}H(f_{t},0))\\
    &\leq 2b[\mathcal{D}_{m}H]_{1,\infty}|p|\cdot|v|-2bv\cdot\mathcal{D}_{m}H(f_{t},0)\nonumber\\
    &\leq b[\mathcal{D}_{m}H]_{1,\infty}(|p|^{2}+|v|^2)-2bv\cdot\mathcal{D}_{m}H(f_{t},0).\nonumber
\end{align}
Collecting all terms and using Young’s inequality we obtain, with $a=\lambda_{B}b$ and $\eta=\gamma+[D]_{1}$,
\begin{align}
    \frac{d}{dt}\int\bigg(\lambda_{B}b|p|^2+2p\cdot v+b|v|^2\bigg)f_{t}(dpdv)&\leq 2bd+(2[A]_{1}+2[D]_{1}b)\int|p|\cdot|v|f_{t}(dpdv)\\
    &+(b\gamma+2\gamma-2\lambda_{B}+2[D]_{1})\int|p|^{2}f_{t}(dpdv)\nonumber\\
    &+(2+\gamma b-2\lambda_{A}b)\int|v|^{2}f_{t}(dpdv)\nonumber\\
    &-2(A(0)+D(0))\cdot\bigg(\int pf_{t}(dpdv)+b\int vf_{t}(dpdv)\bigg)\nonumber\\
    &\leq 2bd+\bigg(\frac{2}{\varepsilon}+\frac{b^{2}\varepsilon}{2[A]^2_{1}}\bigg)|A(0)+D(0)|^{2}\nonumber\\
    &-(2\lambda_{B}-[A]_{1}-2\eta-\varepsilon-\eta b)\int|p|^{2}f_{t}(dpdv)\nonumber\\
    &-\bigg((2\lambda_{A}-\eta )b-2-[A]_{1}-\frac{4[A]^{2}_{1}}{\varepsilon}\bigg)\int|v|^{2}f_{t}(dpdv)\nonumber
\end{align}
for all positive $\varepsilon$. Now, as in the proof of \cref{kineticw2contraction}, with $[A]_{1}$ replaced by $2[A]_{1}$,we get the existence of a positive constant $\eta_{0}$, depending only on $[A]_{1}$, $\lambda_{A}$ and $\lambda_{B}$, such that for all $0\leq\gamma+[D]_{1}<\eta_{0}$, there exist $b$ (and $\varepsilon$) such that $Q_{\lambda_{B}b,b}: (p,v)\mapsto\lambda_{B}b|p|^2+2p\cdot v+b|v|^2$ be a positive quadratic form on $\mathbb{R}^{2d}$ and such that
\begin{align}
    \frac{d}{dt}\int Q_{\lambda_{B}b,b}(p,v)f_{t}(dpdv)&\leq C_{1}-C_{2}\int(|p|^{2}+|v|^{2})f_{t}(dpdv)\\
    &\leq C_{1}-C_{3}\int Q_{\lambda_{B}b,b}(p,v)f_{t}(dpdv)\nonumber
\end{align}
for positive constants $C_{i}$. It follows that
\begin{equation}
    \sup_{t\geq0}\int Q_{\lambda_{B}b,b}(p,v)f_{t}(dpdv)<+\infty
\end{equation}
if initially
\begin{equation}
    \int Q_{\lambda_{B}b,b}(p,v)f_{0}(dpdv)<+\infty,
\end{equation}
that is,
\begin{equation}
    \sup_{t\geq0}\int(|p|^{2}+|v|^{2})f_{t}(dpdv)<+\infty
\end{equation}
if initially $f_{0}$ belongs to $\mathcal{P}_{2}(\mathbb{R}^{d}\times\mathbb{R}^{d})$.\newline

We now turn to the proof of \cref{uniformpropachaoskineticmeanfieldLangevin}. For each $i\in\{1,\ldots,n\}$, the law $f_{t}$ of $(\overline{P}^{(i)}_{t},\overline{V}^{(i)}_{t})$ is the solution to \cref{meanfieldkineticPDE} with $f_{0}$ as initial
datum and the processes $((\overline{P}^{(i)}_{t},\overline{V}^{(i)}_{t}))_{t\geq0}$ and $((P^{(i)}_{t},V^{(i)}_{t}))_{t\geq0}$ are driven by the same Brownian motion. In particular the differences $p^{(i)}_{t}:=P^{(i)}_{t}-\overline{P}^{(i)}_{t}$ and $v^{(i)}_{t}:=V^{(i)}_{t}-\overline{V}^{(i)}_{t}$   
evolve according to
\begin{equation}
        \begin{cases}
            dp^{(i)}_{t}&=v^{(i)}_{t}dt;\\
            dv^{(i)}_{t}&=-\bigg(A(V^{(i)}_{t})-A(\overline{V}^{(i)}_{t})+\lambda_{B}p^{(i)}_{t}+D(P^{(i)}_{t})-D(\overline{P}^{(i)}_{t})+\nabla\frac{\delta H}{\delta m}(\mu_{\mathbf{P}_{t}},P^{(i)}_{t})-\nabla\frac{\delta H}{\delta m}(\mathbb{P}_{P_{t}},\overline{P}^{(i)}_{t})\bigg)dt;\\
            (p^{(i)}_{0},v^{(i)}_{0})&=(0,0).
        \end{cases}
    \end{equation}
As in the proof of \cref{kineticw2contraction}, by the Young inequality and assumptions on $A$ and $D$, for all positive $\varepsilon$, we have
\begin{align}
    \frac{d}{dt}\bigg(\lambda_{B}b|p^{(i)}_{t}|^{2}+2p^{(i)}_{t}\cdot v^{(i)}_{t}+b|v^{(i)}_{t}|^{2}\bigg)&\leq\bigg(\frac{\varepsilon}{2}-2\lambda_{B}+[A]_{1}+2[D]_{1}+[D]_{1}b\bigg)|p^{(i)}_{t}|^{2}\\
    &+\bigg(2+[A]_{1}+\frac{2[A]^{2}_{1}}{\varepsilon}-2\lambda_{A}b+[D]_{1}b\bigg)|v^{(i)}_{t}|^{2}\nonumber\\
    &-2\bigg(p^{(i)}_{t}+bv^{(i)}_{t}\bigg)\cdot\bigg(\nabla\frac{\delta H}{\delta m}(\mu_{\mathbf{P}_{t}},P^{(i)}_{t})-\nabla\frac{\delta H}{\delta m}(\mathbb{P}_{P_{t}},\overline{P}^{(i)}_{t})\bigg).\nonumber
\end{align}
On the other hand, we have
\begin{align}
    &-2\bigg(p^{(i)}_{t}+bv^{(i)}_{t}\bigg)\cdot\bigg(\nabla\frac{\delta H}{\delta m}(\mu_{\mathbf{P}_{t}},P^{(i)}_{t})-\nabla\frac{\delta H}{\delta m}(\mathbb{P}_{P_{t}},\overline{P}^{(i)}_{t})\bigg)\\
    &=-2\bigg(p^{(i)}_{t}+bv^{(i)}_{t}\bigg)\cdot\bigg(\nabla\frac{\delta H}{\delta m}(\mu_{\mathbf{P}_{t}},P^{(i)}_{t})-\nabla\frac{\delta H}{\delta m}(\mu_{\mathbf{P}_{t}},\overline{P}^{(i)}_{t})\bigg)\nonumber\\
    &-2\bigg(p^{(i)}_{t}+bv^{(i)}_{t}\bigg)\cdot\bigg(\nabla\frac{\delta H}{\delta m}(\mu_{\mathbf{P}_{t}},\overline{P}^{(i)}_{t})-\nabla\frac{\delta H}{\delta m}(\mathbb{P}_{P_{t}},\overline{P}^{(i)}_{t})\bigg).\nonumber
\end{align}
By the Cauchy-Schwarz inequality, the Young inequality and assumptions on $H$, we have
\begin{align}
    &-2\bigg(p^{(i)}_{t}+bv^{(i)}_{t}\bigg)\cdot\bigg(\nabla\frac{\delta H}{\delta m}(\mu_{\mathbf{P}_{t}},P^{(i)}_{t})-\nabla\frac{\delta H}{\delta m}(\mathbb{P}_{P_{t}},\overline{P}^{(i)}_{t})\bigg)\\
    &\leq\bigg(2[\mathcal{D}_{m}H]_{1,\infty}+b[\mathcal{D}_{m}H]_{1,\infty}\bigg)|p^{(i)}_{t}|^{2}+b[\mathcal{D}_{m}H]_{1,\infty}|v^{(i)}_{t}|^{2}\nonumber\\
    &+2\bigg(|p^{(i)}_{t}|+b|v^{(i)}_{t}|\bigg)\bigg\|\nabla\frac{\delta H}{\delta m}(\mu_{\mathbf{P}_{t}},\overline{P}^{(i)}_{t})-\nabla\frac{\delta H}{\delta m}(\mathbb{P}_{P_{t}},\overline{P}^{(i)}_{t})\bigg\|\nonumber\\
    &\leq\bigg(2[\mathcal{D}_{m}H]_{1,\infty}+b[\mathcal{D}_{m}H]_{1,\infty}\bigg)|p^{(i)}_{t}|^{2}+b[\mathcal{D}_{m}H]_{1,\infty}|v^{(i)}_{t}|^{2}\nonumber\\
    &+2\|\mathcal{D}^{2}_{m}H\|_{\mathbf{op},\infty}\bigg(|p^{(i)}_{t}|+b|v^{(i)}_{t}|\bigg)\mathcal{W}_{1}(\mu_{\mathbf{P}_{t}},\mathbb{P}_{P_{t}}).\nonumber
\end{align}
By the triangle inequality, we have
\begin{equation}
    \mathcal{W}_{1}(\mu_{\mathbf{P}_{t}},\mathbb{P}_{P_{t}})\leq\mathcal{W}_{1}(\mu_{\mathbf{P}_{t}},\mu_{\mathbf{\overline{P}}_{t}})+\mathcal{W}_{1}(\mu_{\mathbf{\overline{P}}_{t}},\mathbb{P}_{P_{t}}).
\end{equation}
By the Kantorovitch-Rubinstein duality relation, we have
\begin{equation}
    \mathcal{W}_{1}(\mu_{\mathbf{P}_{t}},\mu_{\mathbf{\overline{P}}_{t}})\leq\frac{1}{n}\sum_{j=1}^{n}|p^{(j)}_{t}|.
\end{equation}
Moreover, by symmetry, we have
\begin{equation}
    \mathbb{E}[\mathcal{W}_{1}(\mu_{\mathbf{P}_{t}},\mu_{\mathbf{\overline{P}}_{t}})]\leq\mathbb{E}[|p^{(i)}_{t}|].
\end{equation}
From the Cauchy-Schwarz, Jensen and Young inequalities, we deduce that
\begin{align}
    &2\|\mathcal{D}^{2}_{m}H\|_{\mathbf{op},\infty}\bigg(|p^{(i)}_{t}|+b|v^{(i)}_{t}|\bigg)\mathcal{W}_{1}(\mu_{\mathbf{P}_{t}},\mathbb{P}_{P_{t}})\\
    &\leq\|\mathcal{D}^{2}_{m}H\|_{\mathbf{op},\infty}\bigg(|p^{(i)}_{t}|^{2}+b|v^{(i)}_{t}|^{2}+(b+1)\mathcal{W}_{1}(\mu_{\mathbf{P}_{t}},\mathbb{P}_{P_{t}})^{2}\bigg)\nonumber;\\
    &\mathbb{E}\bigg[2\bigg(|p^{(i)}_{t}|+b|v^{(i)}_{t}|\bigg)\mathcal{W}_{1}(\mu_{\mathbf{P}_{t}},\mathbb{P}_{P_{t}})\bigg]\\
    &\leq(2b+3)\mathbb{E}[|p^{(i)}_{t}|^{2}]+b\mathbb{E}[|v^{(i)}_{t}|^{2}]+(2b+2)\mathbb{E}[\mathcal{W}_{2}(\mu_{\mathbf{\overline{P}}_{t}},\mathbb{P}_{P_{t}})^{2}].\nonumber
\end{align}
By the result of Fournier and Guillin (\cite[Theorem 1]{fournier2013rateconvergencewassersteindistance}), we have
\begin{equation}
    \sup_{t\geq0}\mathbb{E}[\mathcal{W}_{2}(\mu_{\mathbf{\overline{P}}_{t}},\mathbb{P}_{P_{t}})^{2}]\leq C(d)\bigg(\sup_{t\geq0}\int(|p|^{2}+|v|^{2})f_{t}(dpdv)\bigg)\delta_{d}(n).
\end{equation}
Collecting all terms, it follows that there
exists a positive constant $c$ such that for all $\gamma$ and $[D]_{1}$ in $[0,c)$, there exists a constant $M$ such that 
\begin{align}
    \frac{d}{dt}\mathbb{E}\bigg[\bigg(\lambda_{B}b|p^{(i)}_{t}|^{2}+2p^{(i)}_{t}\cdot v^{(i)}_{t}+b|v^{(i)}_{t}|^{2}\bigg)\bigg]&\leq-\bigg(2\lambda_{B}-[A]_{1}-2\eta-\varepsilon-\eta b\bigg)\mathbb{E}[|p^{(i)}_{t}|^{2}]\\
    &-\bigg((2\lambda_{A}-\eta )b-2-[A]_{1}-\frac{4[A]^{2}_{1}}{\varepsilon}\bigg)\mathbb{E}[|v^{(i)}_{t}|^{2}]\nonumber\\
    &+M\bigg(\frac{2}{\varepsilon}+\frac{b^{2}\varepsilon}{2[A]^2_{1}}\bigg)\delta_{d}(n)\nonumber
\end{align}
for all positive $t$, $b$ and $\varepsilon$, where $\eta=\gamma+[D]_{1}$.

Now, as in the proof of \cref{kineticw2contraction}, with $[A]_{1}$ replaced by $2[A]_{1}$,we get the existence of a positive constant $\eta_{0}$, depending only on $[A]_{1}$, $\lambda_{A}$ and $\lambda_{B}$, such that for all $0\leq\gamma+[D]_{1}<\eta_{0}$, there exist $b$ (and $\varepsilon$) such that $Q_{\lambda_{B}b,b}: (p,v)\mapsto\lambda_{B}b|p|^2+2p\cdot v+b|v|^2$ be a positive quadratic form on $\mathbb{R}^{2d}$ and such that
\begin{equation}
    \frac{d}{dt}\mathbb{E}[Q_{\lambda_{B}b,b}(p^{(i)}_{t},v^{(i)}_{t})]\leq-C_{1}\mathbb{E}[|p^{(i)}_{t}|^{2}+|v^{(i)}_{t}|^{2}]+C_{2}\delta_{d}(n)
\end{equation}
for all $t\geq0$ and for positive constants $C_{1}$ and $C_{2}$, also depending on $f_{0}$ through its second moment, but not on $n$. In turn this is bounded by
\begin{equation}
    -C_{3}\mathbb{E}[Q_{\lambda_{B}b,b}(p^{(i)}_{t},v^{(i)}_{t})]+C_{2}\delta_{d}(n),
\end{equation}
so that
\begin{equation}
    \sup_{t\geq0}\mathbb{E}[Q_{\lambda_{B}b,b}(p^{(i)}_{t},v^{(i)}_{t})]\leq C_{4}\delta_{d}(n)
\end{equation}
and finally
\begin{equation}
    \sup_{t\geq0}\mathbb{E}[|p^{(i)}_{t}|^{2}+|v^{(i)}_{t}|^{2}]\leq C\delta_{d}(n)
\end{equation}
where the constant $C$ depends on the parameters of the equation and on the second moment of $f_{0}$, but not on $n$.

\end{proof}

\subsection{Proof related to the examples}
\begin{proof}[Proof of \cref{caracterisationIM}]\label{thm2}
To do this, simply show the two equivalents
$\mathbf{(i)}\Longleftrightarrow\mathbf{(ii)}$ and $\mathbf{(ii)}\Longleftrightarrow\mathbf{(iii)}$.\newline

$\mathbf{(i)}\Longleftrightarrow\mathbf{(ii)}$. Let $\mu_{\infty}$ be such that
\begin{equation}
    \Delta\mu_{\infty}+\nabla\cdot(\mu_{\infty}\mathcal{D}_{m}H(\mu_{\infty},\cdot))=0.
\end{equation}
For any test function $\varphi$, we have
\begin{equation}
    \langle\Delta\mu_{\infty}+\nabla\cdot(\mu_{\infty}\mathcal{D}_{m}H(\mu_{\infty},\cdot)),\varphi\rangle=\int\Delta\varphi d\mu_{\infty}-\int\mathcal{D}_{m}H(\mu_{\infty},\cdot)\cdot\nabla\varphi d\mu_{\infty}.
\end{equation}
Since the Laplacian is a self-adjoint operator, we have
\begin{align}
    \int\Delta\varphi d\mu_{\infty}-\int\mathcal{D}_{m}H(\mu_{\infty},\cdot)\cdot\nabla\varphi d\mu_{\infty}&=-\int\bigg(\nabla\frac{d\mu_{\infty}}{dx}\nabla\varphi+\frac{d\mu_{\infty}}{dx}\mathcal{D}_{m}H(\mu_{\infty},\cdot)\cdot\nabla\varphi\bigg) dx\\
    &=-\int\nabla\bigg(\ln\frac{d\mu_{\infty}}{dx}+\frac{\delta H}{\delta m}(\mu_{\infty},\cdot)\bigg)\cdot\nabla\varphi d\mu_{\infty}.\nonumber
\end{align}
We deduce that
\begin{equation}
    \mathbf{(i)}\Longleftrightarrow\forall\varphi,\quad\int\nabla\bigg(\ln\frac{d\mu_{\infty}}{dx}+\frac{\delta H}{\delta m}(\mu_{\infty},\cdot)\bigg)\cdot\nabla\varphi d\mu_{\infty}\Longleftrightarrow\nabla\bigg(\ln\frac{d\mu_{\infty}}{dx}+\frac{\delta H}{\delta m}(\mu_{\infty},\cdot)\bigg)=0,
\end{equation}
which proves the expected equivalence.\newline

$\mathbf{(ii)}\Longleftrightarrow\mathbf{(iii)}$. As for this last equivalence, we have 
\begin{align}
    \forall (y,z)\in\mathbb{R}^{d}\times\mathbb{R}^{d},\quad&\bigg|\ln\frac{d\mu_{\infty}}{dx}(y)+\frac{\delta H}{\delta m}(\mu_{\infty},y)-\ln\frac{d\mu_{\infty}}{dx}(z)-\frac{\delta H}{\delta m}(\mu_{\infty},z)\bigg|\\
    &\leq \bigg\|\nabla\bigg(\ln\frac{d\mu_{\infty}}{dx}+\frac{\delta H}{\delta m}(\mu_{\infty},\cdot)\bigg)\bigg\|_{\infty}\|y-z\|.\nonumber
\end{align}
We deduce that
\begin{equation}
    \mathbf{(ii)}\Longleftrightarrow \ln\frac{d\mu_{\infty}}{dx}+\frac{\delta H}{\delta m}(\mu_{\infty},\cdot)=C_{\mu_{\infty}}\Longleftrightarrow\mathbf{(iii)}.
\end{equation}

\end{proof}

\begin{proof}[Proof of \cref{proposition1Example1}]\label{proofpropo1example1}
We have
\begin{align}
    \mathcal{D}_{m}H_{VW}(\mu,x)&=\nabla V(x)+\int\nabla W(x-y)\mu(dy);\\
        \mathcal{D}^2_{m}H_{VW}(\mu,x,y)&=-\nabla^2W(x-y).
\end{align}
\end{proof}

\begin{proof}[Proof of \cref{proposition2Example1}]\label{proofpropo2example1}
 We have
\begin{align}
    \mathcal{D}_{m}H_{VW}(\mu,x)&=\nabla V(x)+\int\nabla W(x-y)\mu(dy);\\
        \mathcal{D}^2_{m}H_{VW}(\mu,x,y)&=-\nabla^2W(x-y).
\end{align}   
\end{proof}

\begin{proof}[Proof of \cref{proposition1Example2}]\label{proofpropo1example2}
    The intrinsic derivatives associated with this functional are given by
\begin{align}
   &\mathcal{D}_{m}H^{N}_{VW}(\mu,x)\\
   &=\nabla V(x)\nonumber\\
   &+ \sum_{k=2}^{N}\sum_{j=1}^{k} \int\nabla_{x_{j}}W^{(k)}(x_{1},\ldots,x_{j-1},x,x_{j+1},\ldots,x_{k})\mu^{\otimes k-1}(dx_{1},\ldots,dx_{j-1},dx_{j+1},\ldots,dx_{k})\nonumber\\
   &=\nabla V(x)+\sum_{k=2}^{N}k\int\nabla_{x_{1}}W^{(k)}(x,y)\mu^{\otimes k-1}(dy)\nonumber\\
   \mathcal{D}^2_{m}H^{N}_{VW}(\mu,x,y)&=\sum_{k=2}^{N}k(k-1)\int\nabla^{2}_{x_{1},x_{2}}W^{(k)}(x,y,z)\mu^{\otimes k-2}(dz).
\end{align}
\end{proof}

\begin{proof}[Proof of \cref{proposition2Example2}]\label{proofpropo2example2}
The intrinsic derivatives associated with this functional are given by
\begin{align}
   \mathcal{D}_{m}H^{N}_{VW}(\mu,x)&=\nabla V(x)\\
   &+ \sum_{k=2}^{N}\sum_{j=1}^{k} \int\nabla_{x_{j}}W^{(k)}(x_{1},\ldots,x_{j-1},x,x_{j+1},\ldots,x_{k})\mu^{\otimes k-1}(dx_{1},\ldots,dx_{j-1},dx_{j+1},\ldots,dx_{k})\nonumber\\
   &=\nabla V(x)+\sum_{k=2}^{N}k\int\nabla_{x_{1}}W^{(k)}(x,y)\mu^{\otimes k-1}(dy)\nonumber\\
   \mathcal{D}^2_{m}H^{N}_{VW}(\mu,x,y)&=\sum_{k=2}^{N}k(k-1)\int\nabla^{2}_{x_{1},x_{2}}W^{(k)}(x,y,z)\mu^{\otimes k-2}(dz).
\end{align}
\end{proof}

\begin{proof}[Proof of \cref{proposition1Example3}]\label{proofpropo1example3}
We have
\begin{align}
    \mathcal{D}_{m}H_{\psi}(\mu,x)&=\psi'(\langle\mu,W\rangle)\nabla W(x);\\
        \mathcal{D}^2_{m}H_{VW}(\mu,x,y)&=\psi''(\langle\mu,W\rangle)\nabla W(x)\otimes\nabla W(y).
\end{align}
\end{proof}

\begin{proof}[Proof of \cref{proposition2Example3}]\label{proofpropo2example3}
 We have
\begin{align}
    \mathcal{D}_{m}H_{\psi}(\mu,x)&=\psi'(\langle\mu,W\rangle)\nabla W(x);\\
        \mathcal{D}^2_{m}H_{VW}(\mu,x,y)&=\psi''(\langle\mu,W\rangle)\nabla W(x)\otimes\nabla W(y).
\end{align}   
\end{proof}

\addcontentsline{toc}{section}{References}
\bibliographystyle{unsrt}  
\bibliography{references}

\end{document}